\pdfoutput=1
\documentclass[11pt,a4paper]{article}
\usepackage[T1]{fontenc}
\usepackage{lmodern}
\usepackage[english]{babel}
\usepackage[a4paper,margin=24mm]{geometry}
\usepackage{amsmath,amssymb,amsthm,mathtools,booktabs,array}
\usepackage{microtype,xcolor,enumitem,tikz,caption,placeins,needspace}
\usetikzlibrary{arrows.meta,decorations.markings,positioning}
\makeatletter
\renewcommand*\l@section[2]{%
  \ifnum \c@tocdepth >\z@
    \addpenalty\@secpenalty
    \addvspace{0.12em \@plus\p@}%
    \setlength\@tempdima{1.5em}%
    \begingroup
      \parindent \z@ \rightskip \@pnumwidth
      \parfillskip -\@pnumwidth
      \leavevmode \bfseries
      \advance\leftskip\@tempdima
      \hskip -\leftskip
      #1\nobreak\hfil \nobreak\hb@xt@\@pnumwidth{\hss #2}\par
    \endgroup
  \fi}
\makeatother
\usepackage[colorlinks=true,linkcolor=black,citecolor=black,urlcolor=blue!40!black]{hyperref}
\hypersetup{
  pdftitle={Floquet multiplicities and Krein signatures of elliptic relative equilibria},
  pdfauthor={Xijun Hu; Li Wu},
  pdfsubject={Periodic instability, complex-time monodromy and Floquet multiplicities},
  pdflang={en},
  bookmarksnumbered
}
\allowdisplaybreaks[1]
\setlist{itemsep=3pt,topsep=5pt}
\numberwithin{equation}{section}

\newtheorem{theorem}{Theorem}[section]
\newtheorem{proposition}[theorem]{Proposition}
\newtheorem{lemma}[theorem]{Lemma}
\newtheorem{corollary}[theorem]{Corollary}
\theoremstyle{remark}
\newtheorem{remark}[theorem]{Remark}

\newtheoremstyle{mainresults}
  {10pt}
  {10pt}
  {\itshape}
  {}
  {\bfseries}
  {.}
  {\newline}
  {\thmname{#1}\thmnumber{ #2}\thmnote{ --- #3}}

\theoremstyle{mainresults}
\newtheorem*{introA}{Theorem A}
\newtheorem*{introB}{Theorem B}
\newtheorem*{introCorB}{Corollary B.1}
\newtheorem*{introC}{Theorem C}
\newtheorem*{introCorC}{Corollary C.1}
\newtheorem*{introD}{Theorem D}
\theoremstyle{plain}

\DeclareMathOperator{\rank}{rank}
\DeclareMathOperator{\tr}{tr}
\DeclareMathOperator{\im}{im}
\DeclareMathOperator{\radform}{rad}
\DeclareMathOperator{\diag}{diag}
\DeclareMathOperator{\adj}{adj}
\newcommand{\ii}{\mathrm i}
\newcommand{\dd}{\,\mathrm d}
\newcommand{\RR}{\mathbb R}
\newcommand{\CC}{\mathbb C}
\newcommand{\SSone}{\mathbb S^1}
\newcommand{\vang}{\vartheta}
\newcommand{\q}{q_e}
\newcommand{\ru}{\mathsf u}
\newcommand{\rv}{\mathsf v}
\newcommand{\imag}{\operatorname{Im}}
\newcommand{\norm}[1]{\lVert#1\rVert}

\title{\vspace*{-1.8em}\bfseries Floquet multiplicities and Krein signatures\\
of elliptic relative equilibria}
\author{Xijun Hu\thanks{Partially supported by the National Natural Science
Foundation of China (Grant No.~12521001) and the Taishan Scholars Climbing
Program of Shandong Province (Grant No.~TSPD20240802).}
\and Li Wu\thanks{Partially supported by the National Natural Science
Foundation of China (Grant No.~12171281).}}
\date{}

\begin{document}
\maketitle
\vspace{-2.0em}

\begin{abstract}
We study Floquet multiplicities, Krein signatures, and spectral
instability for elliptic relative equilibria of the planar Newtonian $N$-body
problem. Both the spectral and complex-time arguments stem from the unique
interior zero $r=(-1+\sqrt{1-e^2})/e\in(-1,0)$ of
$q_e(z)=1+\frac{e}{2}(z+z^{-1})$ in the unit disk: a M\"obius angle map sending
$r$ to $0$ makes the twisted Fourier operator tridiagonal and splits its
periodic fiber orthogonally into the constant mode and two conjugate Hardy
half-lines, while a conserved current fixes the Jordan structure at $+1$. In
any even position dimension, removing $(z-1)^{a_1}$ from $\det(zI-M_0)$ leaves
sign $\operatorname{sgn}\det R=(-1)^{N_>(M_0)}$ at $+1$; hence for every
collision-free planar central configuration of positive masses with nonsingular
shape stiffness, the essential multiplier $+1$ is semisimple whenever present,
and an odd shape Morse index forces a positive real multiplier above $1$ for
all $0\le e<1$, including periodic resonances.

For two-position blocks, continuation in $z=e^{\mathrm{i}t}$ factors $M=U_0U_r$
across $0$ and $r$; a Cayley matching pencil and positive transport on $(r,0)$
bound nonperiodic geometric multiplicities by one and separate the Euler and
Lagrange antiperiodic curves. For complex blocks with positive diagonal entries
$a\ne b$ and nonzero coupling, multipliers other than $+1$ have geometric
multiplicity at most one; off $\mathcal T(H)=0$, $\pm1$ are algebraically
simple when present, whereas on $\mathcal T(H)=0$ ($\det H\ne0$) the inner
Hardy kernel forces semisimple $(1,1)$ Krein-double periodic roots. For
centered regular $(1+n)$-gons ($0<e<1$), periodic roots are simple in mode $1$
($n=4$, $n\ge7$), mode $2$ ($n\ge5$), and modes $4\le\ell<n/2$, while mode $3$
realizes infinitely many $(1,1)$ double roots $e_j\uparrow1$ for $n\ge100$.
\end{abstract}


\setcounter{tocdepth}{1}
\vspace{-0.6em}
\begingroup
\footnotesize
\renewcommand{\contentsname}{\normalsize Contents}
\tableofcontents
\endgroup
\medskip

\section{Introduction and main results}\label{sec:intro}

An elliptic relative equilibrium of the planar Newtonian $N$-body problem
rotates and scales a central configuration along a Kepler ellipse of
eccentricity $e\in[0,1)$. Passing to rotating--pulsating coordinates and
removing translations and the homographic Kepler motions leaves a real
symplectic return $M_0$ on the essential shape variations. At the periodic multiplier $+1$, the question is whether a
resonance can change the parity that forces a positive real multiplier above $1$.
For $0<e<1$, the continued Kepler weight
$\q(z)=1+\frac{e}{2}(z+z^{-1})$ has a unique zero
$r=(-1+\sqrt{1-e^2})/e$ inside the unit disk. The disk automorphism sending
$r$ to the origin puts the variational operator in tridiagonal Fourier form.
At the periodic twist both couplings to the constant mode vanish. The
resulting orthogonal splitting into the constant mode and two Hardy
half-lines also holds at $e=0$, in arbitrary even position dimension.
Together with the conserved current, it determines the Jordan structure
at $+1$ and yields the parity law for nonsingular shape stiffness.

For two-position blocks, a second question concerns the dimension of a
Floquet eigenspace and, in the Euler and Lagrange families, the separation
of the antiperiodic degeneracy curves. For $0<e<1$, continuation around
the two interior singularities $0$ and $r$ factors the return as
$M=U_0U_r$. Transport along $(r,0)$ compares the local solutions at the two
ends. Its sign information enters the matching pencil and bounds the
nullity at nonperiodic twists. The two arguments meet at a periodic
resonance on $\mathcal T(H)=0$ for a coupled complex block. There the
common local fixed space $\ker A\cap\ker B$ coincides with the inner
Hardy kernel, and the periodic multiplier is semisimple double with
Krein signature $(1,1)$.

\subsection{The periodic return on planar shapes}\label{subsec:headline}
Let $J^T=-J$, $J^2=-I$, and $R=R^T$ on an even-dimensional real position space,
put $\q(t)=1+e\cos t$, and let $M_0$ be the $2\pi$-return on position and
velocity of
\[
 w''+2Jw'=\q^{-1}Rw .
\]
At a collision-free planar Newtonian central configuration $a_0$ of positive
masses, normalized by $\nabla_m U(a_0)=-\lambda_c a_0$, the essential shape
stiffness is $R=I+\lambda_c^{-1}D\nabla_m U(a_0)$ restricted to the
mass-orthogonal complement of translations, $a_0$, and $Ja_0$. Its Morse index
$n_-(R)$ counts the negative directions of the Newtonian potential on the shape
sphere transverse to rigid rotations.

\begin{introA}[Periodic parity and semisimplicity on planar shapes]
\phantomsection\label{main:parity}
Let $R$ be nonsingular and $0\le e<1$. If $a_1$ denotes the algebraic
multiplicity of $+1$ in $M_0$ and $p_*(z)=\det(zI-M_0)/(z-1)^{a_1}$, then
\[
 \operatorname{sgn}p_*(1)=\operatorname{sgn}\det R=(-1)^{N_>(M_0)},
\]
where $N_>(M_0)$ counts the positive real multipliers greater than one with
algebraic multiplicity. For every collision-free planar Newtonian central
configuration of positive masses with invertible shape stiffness, the essential
multiplier $+1$ is semisimple whenever it occurs; consequently, an odd shape
Morse index forces a positive real multiplier greater than one at every
eccentricity $0\le e<1$, including periodic resonances.
\end{introA}

For circular relative equilibria Roberts derived a necessary condition for spectral stability from the linearized spectrum \cite{Rob99}. Odd-index instability for circular relative equilibria ($e=0$) was proved by Hu
and Sun \cite{HS} and further analyzed by Asselle, Portaluri, and Wu
\cite{APW}. At positive eccentricity, Mart\'inez, Sam\`a, and Sim\'o \cite{MSS}
established a Fourier pairing for a real four-state family, Zhou and Long
\cite{ZL} proved the persistence of an expanding real pair across periodic
resonances in the elliptic Euler family, and Chu \cite{Chu} obtained the
full-eccentricity expanding-pair result in the indefinite region of a real
four-state restricted problem via a weighted Fourier pairing at $+1$.
\hyperref[main:parity]{Theorem~A} extends the periodic parity law to arbitrary
even position dimension and establishes universal semisimplicity at $+1$ for
the full essential planar Newtonian return.

The proof in Sections~\ref{sec:graph} and~\ref{sec:periodic-general} works
through the M\"obius-adapted tridiagonal operator for an arbitrary unit twist
$\lambda=e^{2\pi\ii\theta}$. Eliminating its coercive high-frequency tail
produces a finite Hermitian Schur complement $S(\theta,x)$ whose determinant
matches $\det(\lambda I-M_x)$ up to a nonvanishing real-analytic factor,
derived via the relative Hill trace formula in the setting of Hu, Ou, and Wang
\cite{HOW}. Thus the vanishing order of $\det S$ in $\theta$ gives the
algebraic multiplicity of $\lambda$, its order in a scalar stiffness shift $x$
gives the geometric multiplicity, and $\partial_\theta S|_{\ker S}$ recovers
the conserved Hermitian current (Propositions~\ref{prop:graph-twist}
and~\ref{prop:graph-det}), realizing the self-adjoint pencil formulation of the
Krein signature discussed by Koll\'ar and Miller \cite{KM}. At $\theta=0$, the
constant mode decouples completely from the two conjugate Hardy half-lines,
identifying the exact Jordan blocks at $+1$ and reducing the real periodic sign
to $\operatorname{sgn}\det R$.

\subsection{Two-position equations and complex-time factorization}\label{subsec:setup}

Let $t$ denote true anomaly, let $e\in[0,1)$ be the eccentricity, and set
$\q(t)=1+e\cos t$. In our rotating-frame conventions, the real and complex
two-component variational equations take the forms
\begin{align}
 L_{D,e}u&=-u''-2Ju'+\q^{-1}Du=0,
 &J&=\begin{pmatrix}0&-1\\1&0\end{pmatrix},\quad D=D^T\in\RR^{2\times2},
 \label{eq:real-op}\\
 \mathcal L_{H,e}w&=-w''+2\ii Zw'+\q^{-1}Hw=0,
 &Z&=\diag(1,-1),\quad H=\begin{pmatrix}a&c\\c&b\end{pmatrix}\in\RR^{2\times2},
 \label{eq:complex-op}
\end{align}
where primes denote differentiation with respect to $t$. Denote by $M$ the
$2\pi$-period map acting on the four-dimensional complex state space (using the
complexification of the real state space for \eqref{eq:real-op}), and define
the geometric Floquet multiplicity at $\lambda\in\CC$ by
\[
 \nu_\lambda(M)=\dim_{\CC}\ker(M-\lambda I).
\]
Whenever $|\lambda|=1$, evaluation of initial data identifies
$\ker(M-\lambda I)$ with the kernel of $L_{D,e}$ or $\mathcal L_{H,e}$ on
$L^2([0,2\pi],\CC^2)$ with operator domain
\[
 \{f\in H^2([0,2\pi],\CC^2):f(2\pi)=\lambda f(0),\ f'(2\pi)=\lambda f'(0)\},
\]
whose associated quadratic form has domain
$\{f\in H^1([0,2\pi],\CC^2):f(2\pi)=\lambda f(0)\}$. In polygonal
applications, $\nu_\lambda(M)$ always refers to one specified block; at
$\lambda=\pm 1$, complexification of a real block preserves the numerical
dimension of its kernel.

\begin{figure}[t]
\centering
\begin{tikzpicture}[scale=2.35,>=Stealth,
  midarrow/.style={postaction={decorate},
    decoration={markings,mark=at position #1 with {\arrow{Stealth[length=5pt]}}}}]
  \draw[gray!55,thin,->] (-2.15,0) -- (1.42,0) node[right,black,font=\small] {$\operatorname{Re} z$};
  \draw[gray!55,thin,->] (0,-1.22) -- (0,1.22) node[above,black,font=\small] {$\operatorname{Im} z$};

  \draw[thick,black,midarrow=0.13,midarrow=0.63] (0,0) circle (1);
  \node[font=\small] at (0.94,0.88) {$|z|=1$};
  \node[font=\footnotesize,gray!30!black] at (0.98,-0.88) {$M\sim U_0U_r$};

  \draw[very thick,teal!70!black] (-0.62,0) -- (0,0);

  \draw[blue!65!black,thick,midarrow=0.52]
    (-0.31,0.032) -- (-0.483,0.032)
    arc[start angle=13,end angle=347,radius=0.14]
    -- (-0.31,-0.032);
  \node[blue!65!black,font=\small] at (-0.64,0.26) {$U_r$};

  \draw[red!65!black,thick,midarrow=0.52]
    (-0.31,-0.032) -- (-0.137,-0.032)
    arc[start angle=193,end angle=527,radius=0.14]
    -- (-0.31,0.032);
  \node[red!65!black,font=\small] at (0.14,0.26) {$U_0$};

  \draw[dashed,gray!50!black,thin]
    (-0.31,0.015) .. controls (-0.12,0.68) and (0.52,0.61) .. (1,0);

  \fill[black] (-1.62,0) circle (1.2pt);
  \node[font=\small] at (-1.62,-0.16) {$r^{-1}$};
  \fill[blue!65!black] (-0.62,0) circle (1.3pt);
  \node[blue!65!black,font=\small] at (-0.77,-0.22) {$r$};
  \fill[red!65!black] (0,0) circle (1.3pt);
  \node[red!65!black,font=\small] at (0.17,-0.22) {$0$};
  \fill[teal!70!black] (-0.31,0) circle (1.2pt);
  \node[teal!70!black,font=\small] at (-0.31,-0.16) {$z_b$};
  \fill[black] (1,0) circle (1.2pt);
  \node[font=\small] at (1.13,-0.15) {$1$};
\end{tikzpicture}
\caption{Interior regular singularities $0$ and $r=(-1+\sqrt{1-e^2})/e$ of
\eqref{eq:real-op}--\eqref{eq:complex-op} in the complex time plane
$z=e^{\ii t}$, together with the exterior pole $r^{-1}<-1$, the real
connecting segment $(r,0)$, the common base point $z_b\in(r,0)$, and the two
positively oriented local loops. The ordered product represents the physical
return up to change of base fiber.}
\label{fig:z-plane}
\end{figure}
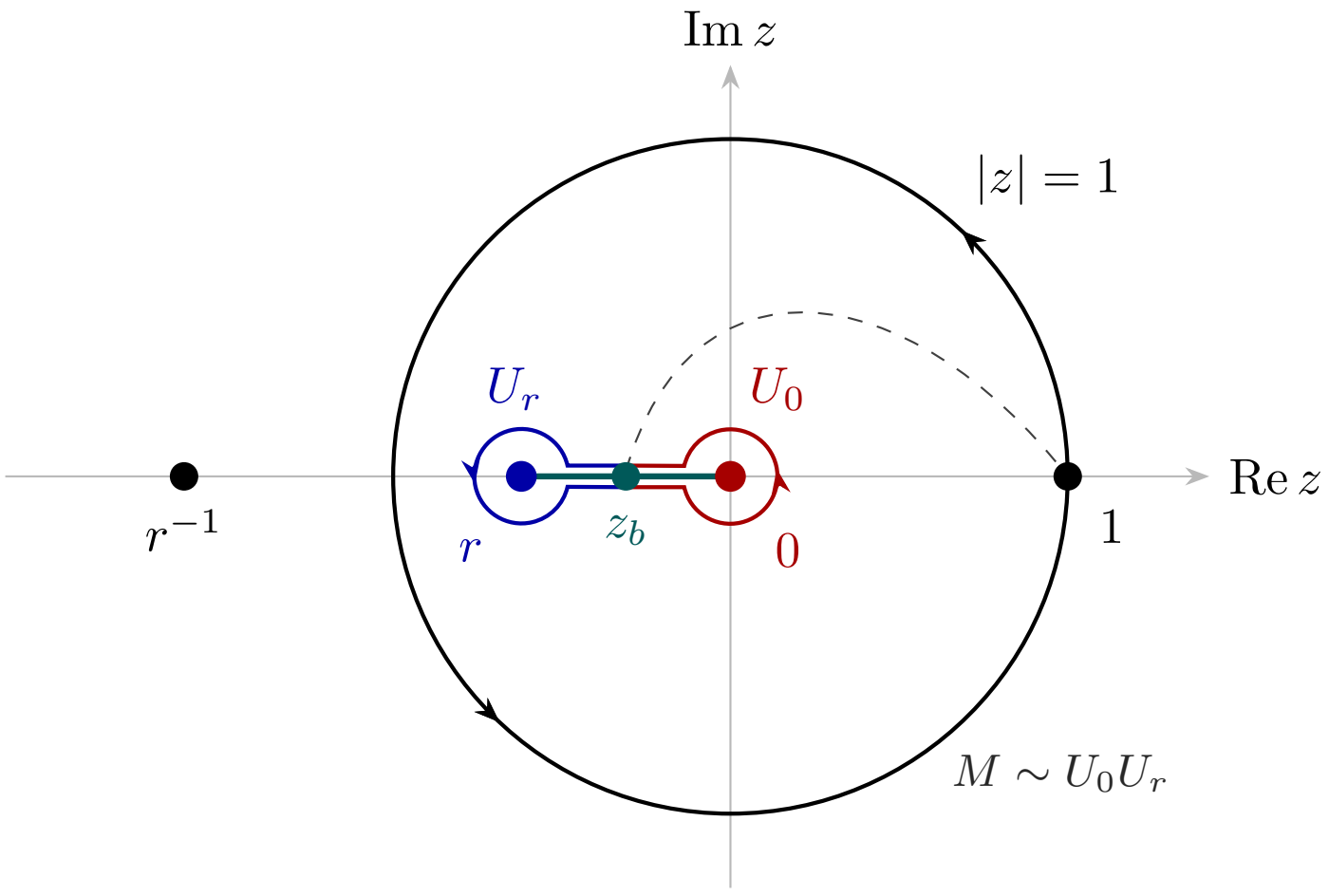

For $0<e<1$, introduce the complex time variable $z=e^{\ii t}$ together with
the Euler derivative $\vang=z\,\mathrm d/\mathrm dz$ and the parameters
\[
 \eta=\sqrt{1-e^2},\qquad r=\frac{-1+\eta}{e}\in(-1,0).
\]
The analytic continuation $\q(z)=1+e(z+z^{-1})/2$ vanishes only at $r$ and at
$r^{-1}<-1$. As $e\uparrow1$ these two zeros approach $-1$ from opposite sides
of the unit circle. The identity $r=-e/(1+\eta)$ gives $r\to0$ as $e\downarrow0$.
For $0<e<1$, $0$ and $r$ are the only possible singularities of the
continued equations inside the disk, and both are of regular-singular type
when present. Fix a common base point $z_b\in(r,0)$ and
let $U_0,U_r$ be the local monodromies along positively oriented loops based at
$z_b$ that approach $0$ and $r$ along $(r,0)$, as shown in
Figure~\ref{fig:z-plane}. Because the positively oriented unit circle is
homotopic in $\CC\setminus\{0,r,r^{-1}\}$ to the concatenation of these two
based loops, identifying the physical fiber at $z=1$ with the fiber at $z_b$
via a compatible stem gives the common-base-point factorization
\begin{equation}\label{eq:product}
 M=U_0U_r.
\end{equation}
Reversing the order of the factors produces a similar matrix. Moreover, in the
$z$-frames introduced in Sections~\ref{sec:real} and~\ref{sec:complex}, both
differential equations have real coefficients along the connecting segment
$(r,0)$.

In both \eqref{eq:real-op} and \eqref{eq:complex-op}, the based local factors
are unipotent, so their Cayley generators $A=-2\ii(U_0-I)(U_0+I)^{-1}$ and
$B=-2\ii(U_r-I)(U_r+I)^{-1}$ satisfy $A^3=0$ and $B^2=0$ (with $A^2=0$ in the
real system). In Section~\ref{sec:theory}, Theorem~\ref{thm:cayley} converts
the eigenvalue equation $(M-\lambda I)v=0$ into the linear matching pencil
$\mathcal C(\lambda)=(1-\lambda)(I-AB/4)+\frac{\ii}{2}(1+\lambda)(A+B)$, whose
determinant agrees with both $\det(M-\lambda I)$ and the finite Schur
determinant of \eqref{eq:same-characteristic} up to a nonvanishing factor. In
four dimensions, any nonperiodic multiplicity $\nu_{\lambda_0}(M)\ge 2$ forces
the characteristic polynomial to be a perfect reciprocal square and hence annihilates
the trace difference $\Delta(M)=\tr M-\tr M^{-1}=-\ii\tr(A^2B)$; when $A^2=0$,
the same obstruction forces $\dot\Delta=0$ along any determinant-preserving
variation.

\subsection{Multiplicity, curve separation, and the polygonal trichotomy}\label{subsec:results}

For the real two-component equation \eqref{eq:real-op}, the local generators
$A,B$ lie in $\mathfrak{sp}(4,\RR)$ and satisfy $A^2=B^2=0$, so $\Delta(M)$
vanishes identically and simplicity is instead extracted from the real
symmetric pencil induced by $\mathcal C(\lambda)$ together with positive
transport along $(r,0)$.

\needspace{5\baselineskip}
\begin{introB}[Real symmetric stiffness]
\phantomsection\label{main:real}
Fix $0<e<1$ in \eqref{eq:real-op}.
\begin{enumerate}
\item If $D$ is indefinite, has rank at most one, or is positive definite and
nonscalar, then $\nu_\lambda(M)\le 1$ for every $\lambda\neq 1$. Moreover, for
indefinite $D$, all Floquet multipliers lie in $\RR\cup\SSone$.
\item The antiperiodic kernel vanishes if $D=\alpha I$ with
$\alpha\ge\min\{1,3(1+e)/4\}$, or if $\norm{D}_{\rm op}\le 3(1-e)/4$.
\end{enumerate}
In particular, $\tr D\ge 2$ implies $\nu_{-1}(M)\le 1$.
\end{introB}

For the elliptic Lagrange triangle, Roberts reduced the essential dynamics to four dimensions and traced the period-doubling boundary \cite{Rob02}. For planar three-body relative equilibria, Hu, Long, and Sun \cite{HLS}
established the antiperiodic counts and analytic branches for the equilateral
Lagrange family, Zhou and Long \cite{ZL} enumerated the extended Euler curves,
and Hu, Wu, and Zhu \cite{HWZ} recently proved strict antiperiodic separation
at positive eccentricity for both families. Combining
\hyperref[main:real]{Theorem~B} with a self-contained derivation of the
three-body shape stiffnesses and an index count from the Kepler endpoint
provides an independent complex-time proof of these separations.

\begin{introCorB}[Euler and Lagrange antiperiodic curves]
\phantomsection\label{main:real-curves}
For every $0<e<1$, the two branches of each antiperiodic pair in the extended
Euler family ($\delta>0$) are strictly ordered, and the Lagrange family
($\beta\in[0,9]$) admits exactly two antiperiodic parameter values, which are
distinct for $0<e<1$ and coalesce at $\beta=3/4$ for $e=0$. At every
positive-eccentricity degeneracy point, the antiperiodic kernel is
one-dimensional and each degeneracy curve is real-analytic in $e\in(0,1)$.
\end{introCorB}

For the single complex block \eqref{eq:complex-op}, the physical conserved
form is Hermitian and the origin generator satisfies $A^3=0$, with $A^2$
controlled by the quadratic stiffness polynomial
\begin{equation}\label{eq:T}
 \mathcal T(H)=3(a-b)^2+6(a+b)-4c^2-9.
\end{equation}
For $a,b>0$, $c(a-b)\neq 0$, and $0<e<1$, pairing with the transpose-dual
equation along $(r,0)$ gives $\Delta(M)\neq 0$ whenever $\mathcal T(H)\neq 0$,
while a transverse stiffness shift $H+\varepsilon I$ yields the variational
obstruction of Theorem~\ref{thm:cayley}(2b) along $\mathcal T(H)=0$.

\begin{introC}[Nonperiodic simplicity for a coupled complex block]
\phantomsection\label{main:complex}
Suppose $a,b>0$, $a\neq b$, and $c\neq 0$ in \eqref{eq:complex-op}. Then for
every $0<e<1$ and every $\lambda\neq 1$, one has $\nu_\lambda(M)\le 1$.
Moreover, whenever $\mathcal T(H)\neq 0$, each of the multipliers $+1$ and
$-1$, if present, has algebraic multiplicity one.
\end{introC}

Hu, Long, and Ou \cite{HLO} decomposed the essential linearized system of a
regular $(1+n)$-gon into real Hamiltonian blocks. Working in circular
coordinates isolates one complex four-state block \eqref{eq:complex-op} at a
time (together with a single real terminal block for even $n$), so
\hyperref[main:complex]{Theorem~C} and \hyperref[main:real]{Theorem~B} yield
uniform nonperiodic and antiperiodic simplicity bounds across all physical
modes.

\begin{introCorC}[Physical modes of a regular $(1+n)$-gon]
\phantomsection\label{main:physical}
Place $n\ge 3$ unit masses at the vertices of a regular polygon and a mass
$m>0$ at its center. For $0<e<1$:
\begin{enumerate}
\item The first complex mode satisfies $\nu_\lambda(M)\le 1$ for all
$\lambda\neq 1$; for $n\ge 7$, its multipliers at $+1$ and $-1$, whenever
present, are algebraically simple.
\item In each higher complex mode $2\le\ell<n/2$, one has
$\nu_\lambda(M)\le 1$ for every $\lambda\neq 1$ at every coupled mass
($m\neq 2Q_\ell$), and $\nu_{-1}(M)\le 1$ at every mass $m>0$.
\item For even $n$, the real four-dimensional terminal mode $\ell=n/2$
satisfies $\nu_{-1}(M)\le 1$ for every $m>0$.
\end{enumerate}
\end{introCorC}

For the coupled unequal-diagonal family of \hyperref[main:complex]{Theorem~C},
a two-dimensional periodic eigenspace can arise only on the resonant surface
$\mathcal T(H)=0$. There the complex-time fixed space $E=\ker A\cap\ker B$
consists of holomorphic solutions vanishing at $r$ and therefore coincides with
the inner Hardy kernel of Section~\ref{sec:periodic-general}. Time reflection
supplies the conjugate outer kernel, Proposition~\ref{thm:current-split} fixes
the $(1,1)$ Krein signature and semisimplicity, and Riesz perturbation at fixed
eccentricity produces opposite mass--phase crossing slopes
(Corollary~\ref{cor:mass-phase}). On the polygonal side, Hu and Ou \cite{HO}
proved divergence of the periodic Morse index of the full essential system
below the collision threshold $-1/8$; our shifted determinant identity on
$\mathcal T(H)=0$ localizes this phenomenon precisely to the third mode, whereas
a general two-matrix positive-sector coercivity criterion
(Lemma~\ref{lem:square-hardy}) rules out periodic roots at the exceptional mass
of the square ($n=4$) and full-domain coercivity excludes them in all higher
modes $\ell\ge 4$.

\begin{introD}[Periodic eigenspaces on $\mathcal T=0$ and the polygonal trichotomy]
\phantomsection\label{main:resonant}
Consider the complex block \eqref{eq:complex-op} and the physical polygonal
normalization of \hyperref[main:physical]{Corollary~C.1}.
\begin{enumerate}
\item Suppose $a,b>0$, $c\neq 0$, $\mathcal T(H)=0$, and $\det H\neq 0$ in
\eqref{eq:complex-op}. Then for every $0<e<1$, the periodic multiplicity
satisfies $\nu_1(M)\in\{0,2\}$. Whenever $\nu_1(M)=2$, $+1$ is a semisimple
eigenvalue of algebraic multiplicity two, and the conserved physical Hermitian
form of Proposition~\ref{prop:current} has signature $(1,1)$ on $\ker(M-I)$.
\item In the first complex mode of the square ($n=4$), one has
$\nu_\lambda(M)\le 1$ for every $m>0$, $0<e<1$, and $\lambda\in\CC$. In
particular, at the unique positive zero $m_*=(839+395\sqrt2)/496$ of
$\mathcal T(H)$, the periodic operator is invertible for all $0\le e<1$.
\item In each higher complex mode $2\le\ell<n/2$, the condition
$\mathcal T(H(m))=0$ determines a unique candidate mass
$m_{\mathcal T}=m_{\mathcal T}(n,\ell)\in\RR$. Whenever $m_{\mathcal T}>0$ and
$\det H(m_{\mathcal T})\neq 0$, every periodic root at $m_{\mathcal T}$ with
$0<e<1$ has a semisimple multiplier $+1$ of algebraic multiplicity two with
Hermitian signature $(1,1)$, and at fixed $e$ the two real mass curves for the
phase $\lambda=e^{2\pi\ii\theta}$ cross with opposite slopes at $\theta=0$. If
moreover $\det(H(m_{\mathcal T})+I/8)<0$, then there exists a sequence of
distinct periodic roots $e_j\uparrow 1$ at $m=m_{\mathcal T}$. Specifically:
\begin{enumerate}
\item For $\ell=2$ and every $n\ge 5$, one has $m_{\mathcal T}<0$; hence every
periodic multiplier at every $m>0$ and $0<e<1$ is algebraically simple.
\item For $\ell=3$ and every $n\ge 100$, one has $m_{\mathcal T}>0$ and
$\det(H(m_{\mathcal T})+I/8)<0$, producing infinitely many eccentricities
$e_j\uparrow 1$ at $m=m_{\mathcal T}$ with semisimple periodic multiplicity
two and Hermitian signature $(1,1)$.
\item For every $4\le\ell<n/2$, one has $m_{\mathcal T}>0$ and
$H(m_{\mathcal T})>0$; at every $m>0$ with $H(m)>0$ the full periodic form is
uniformly coercive in $0\le e<1$ (Proposition~\ref{prop:fixed-modes}). In
particular, the periodic operator at $m_{\mathcal T}$ is invertible for all
$0\le e<1$, and every periodic multiplier at $m>0$, $0<e<1$ is algebraically
simple.
\end{enumerate}
\end{enumerate}
\end{introD}

\subsection{Organization of the paper}\label{subsec:org}

Section~\ref{sec:graph} constructs the M\"obius-adapted tridiagonal Fourier
operator, its finite Hermitian Schur complement, and the relative Hill formula.
Section~\ref{sec:periodic-general} establishes the orthogonal Hardy splitting
at $\theta=0$, the Jordan and current classification, and
\hyperref[main:parity]{Theorem~A}. For reversible systems,
Proposition~\ref{prop:reversible-half} expresses the full characteristic
polynomial and the two endpoint kernels in terms of a forward half-period map.
Section~\ref{sec:theory} develops the
two-factor Cayley matching pencil (Theorem~\ref{thm:cayley}), whose connection
data are evaluated for real and complex two-position systems in
Sections~\ref{sec:real} and~\ref{sec:complex} to prove
\hyperref[main:real]{Theorem~B}, \hyperref[main:real-curves]{Corollary~B.1},
and \hyperref[main:complex]{Theorem~C}. Sections~\ref{sec:polygon}
and~\ref{sec:periodic} treat the regular $(1+n)$-gons and prove
\hyperref[main:physical]{Corollary~C.1} and
\hyperref[main:resonant]{Theorem~D}. Appendix~\ref{app:gauge} proves
convergence of the local normal forms and evaluates the order-two resonance
coefficients; Appendix~\ref{app:hardy} proves the uniform full-domain and
two-matrix positive-sector coercivity criteria together with the circular
spectrum on $\mathcal T=0$; and Appendix~\ref{app:finite-n} collects the
self-contained rational estimates for finite polygons.

\section{Tridiagonal Fourier representation and relative Hill formula}\label{sec:graph}
Let the position space be $\CC^s$ equipped with its Hermitian inner product,
antilinear in the first argument. For a unitary twist
$\lambda=\exp(2\pi\ii\theta)$ with real lift $\theta$, a real scalar stiffness
shift $x$, and matrices $Z=Z^*=Z^{-1}$ and $H=H^*$, consider the boundary-value
problem
\begin{equation}\label{eq:whole-shape}
 w''-2\ii Zw'=\q^{-1}(H+xI)w,\qquad
 (w,w')(2\pi)=\lambda(w,w')(0).
\end{equation}
We use the normalized circle measures $\dd t/(2\pi)$ and $\dd\phi/(2\pi)$ in
all quadratic forms of this section and Section~\ref{sec:periodic-general},
while the physical Green kernel below is taken relative to Lebesgue measure
$\dd t$. Let $M_x$ denote the $2\pi$-return of \eqref{eq:whole-shape} in
position and velocity coordinates. For any two solutions $v,w$ at real $x$,
differentiation shows that the physical current
\begin{equation}\label{eq:common-current}
 [v,w]_{\mathcal J}=-\ii\big(v^*w'-v'^*w-2\ii v^*Zw\big)
\end{equation}
is constant in $t$. In the state coordinates $(w,-\ii w')$, its matrix is
$\bigl(\begin{smallmatrix}-2Z&I\\I&0\end{smallmatrix}\bigr)$, which is
nondegenerate, so $M_x$ preserves the pairing \eqref{eq:common-current}.

\subsection{Eliminating the high Fourier modes}
Set $\eta=\sqrt{1-e^2}$ and $r=-e/(1+\eta)$, including $r=0$ at the circular
orbit. The M\"obius change of angle
\[
 e^{\ii\phi}=\frac{e^{\ii t}-r}{1-r e^{\ii t}},\qquad
 \frac{\dd\phi}{\dd t}=\frac{\eta}{\q(t)}=\frac{1-e\cos\phi}{\eta}
\]
preserves orientation and advances $\phi$ by $2\pi$ over one physical period.
Writing $w(t)=e^{\ii\theta\phi(t)}u(\phi(t))$ with $u$ $2\pi$-periodic and
setting $D_\theta=-\ii\partial_\phi+\theta$, a change of variable in the
twisted $H^1$ quadratic form of $-\partial_t^2+2\ii Z\partial_t+(H+xI)/\q$
followed by multiplication by $\eta$ yields the quadratic form of the
transformed operator
\begin{equation}\label{eq:physical-jacobi}
 T(\theta,x)=D_\theta(1-e\cos\phi)D_\theta-2\eta ZD_\theta+H+xI,
\end{equation}
with periodic form domain $H^1_{\rm per}(\dd\phi/(2\pi),\CC^s)$. The change of variables preserves both the Morse index and the nullity. In the Fourier basis of
$L^2(\dd\phi/(2\pi),\CC^s)$, the shifted frequencies $p_k=k+\theta$ produce
diagonal blocks $p_k^2 I-2\eta p_k Z+H+xI$ and nearest-neighbor blocks
$-ep_kp_{k+1}I/2$, so $T(\theta,x)$ is block tridiagonal.

\begin{proposition}[Finite Schur reduction, inertia, and the phase derivative]\label{prop:graph-twist}
Fix $(e,\theta,x)$ with $0\le e<1$. In a neighborhood of this point, let $P$
project onto a sufficiently large fixed interval of Fourier modes, set
$Q=I-P$, and define
\begin{equation}\label{eq:finite-graph}
 L_{\rm tail}=QTQ|_{Q},\qquad
 A_{\rm c}=PTP|_P,\qquad C=QTP|_P,\qquad
 B_{\rm c}=C^*L_{\rm tail}^{-1}C,\qquad S=A_{\rm c}-B_{\rm c}.
\end{equation}
Then $L_{\rm tail}$ is strictly positive and $S$ is a finite Hermitian matrix
whose negative index coincides with the Morse index of the twisted physical
form and whose kernel maps isomorphically onto $\ker(M_x-\lambda I)$. With
$Xv=v-L_{\rm tail}^{-1}Cv$ and $w_v=e^{\ii\theta\phi}Xv$ in physical time, one
has
\begin{equation}\label{eq:graph-current}
 \partial_x S=X^*X>0,\qquad
 v^*(\partial_\theta S)u=\eta[w_v,w_u]_{\mathcal J}
 \quad(v,u\in\ker S).
\end{equation}
In particular, $\lambda$ is semisimple if and only if the restriction of
$\partial_\theta S$ to $\ker S$ is nondegenerate.
\end{proposition}

\begin{proof}
Choose $P$ to contain all frequencies in a large interval about $-\theta$ at
the base point and hold $P$ fixed locally. Whenever $|p_k|\ge b>0$ on the
tail, the tail quadratic form is bounded below by
\[
 \big(1-e-2\eta/b-(\|H\|+|x|)/b^2\big)\|D_\theta u\|^2,
\]
which is uniformly positive in a local neighborhood once $b$ is taken large
enough. Consequently, $L_{\rm tail}$ has a bounded inverse from the dual of
its $H^1$ form domain back into that domain, and by a Neumann series both the
operator coefficients and $L_{\rm tail}^{-1}$ depend real-analytically on the
local parameters. For $v\in P$ and $y\in QH^1$, completing the square yields
\[
 \langle v+y,T(v+y)\rangle
 =v^*Sv+\|L_{\rm tail}^{1/2}(y+L_{\rm tail}^{-1}Cv)\|^2,
\]
which matches the negative indices and sends each $v\in\ker S$ to the unique
full kernel vector $Xv$ projecting to $v$; integration by parts then recovers
the $H^2$ solution and its derivative boundary condition.

Since $S=X^*TX$, any derivative of $X$ takes values in $Q$ while $TX$ takes
values in $P$, so only the derivative of $T$ contributes; thus $\partial_x T=I$
yields the first identity in \eqref{eq:graph-current}. For the phase
derivative, hold $Xv$ and $Xu$ fixed and differentiate the physical lifts
$e^{\ii\theta\phi}Xv$ and $e^{\ii\theta\phi}Xu$: differentiating the two
derivative terms in the physical form produces the integrand
$\phi'(t)[w_v,w_u]_{\mathcal J}$, where $[w_v,w_u]_{\mathcal J}$ is constant
on solutions and $\phi'$ has mean value one over $[0,2\pi]$, while the factor
$\eta$ comes from \eqref{eq:physical-jacobi}. Finally, because $|\lambda|=1$
and $M_x$ preserves the nondegenerate form $[\,\cdot,\cdot\,]_{\mathcal J}$
(so that $M_x^\dagger=M_x^{-1}$), taking the $\mathcal J$-adjoint yields
\[
 \operatorname{ran}(M_x-\lambda I)=\ker(M_x-\lambda I)^{\perp_{\mathcal J}}.
\]
The intersection of this range with $\ker(M_x-\lambda I)$ consists precisely of
those eigenvectors admitting a Jordan predecessor, which establishes the
semisimplicity criterion.
\end{proof}

At a simple root, the implicit function theorem gives the phase derivative of
the scalar stiffness shift as
\[
 \frac{\dd x}{\dd\theta}=-\frac{\eta[w_v,w_v]_{\mathcal J}}{\|Xv\|^2},
 \qquad 0\ne v\in\ker S.
\]
Enlarging the Fourier interval $P$ adds only positive directions to the square
completion, leaving both the Morse index and the conserved current on physical
solutions independent of the cutoff. For a real physical equation, complex
conjugation exchanges the twists $\lambda$ and $\bar\lambda$; at
$\lambda=\pm 1$, choosing a frequency interval invariant under
$k\mapsto -k-2\theta$ makes the reduction commute with conjugation, and
complexifying its restriction to the real subspace preserves the inertia,
whereas for a complex modal block $S$ is already Hermitian on that block alone.

\subsection{Relative Hill formula and the local determinant}

\begin{proposition}[Relative Hill formula for the reduced determinant]\label{prop:graph-det}
Near each real $(\theta_0,x_0)$ with $e<1$ fixed, there exists a nowhere-zero
real-analytic function $g(\theta,x)$ such that
\begin{equation}\label{eq:same-characteristic}
 \det(e^{2\pi\ii\theta}I-M_x)=g(\theta,x)\det S(\theta,x).
\end{equation}
In particular, the vanishing order of $\det S$ in $\theta$ at fixed stiffness
equals the algebraic multiplicity of the multiplier, whereas its order in $x$
at fixed phase equals the geometric multiplicity.
\end{proposition}
\begin{proof}
Fix a real shift $x_*>x_0$ large enough that $T(\theta,x_*)>0$ locally in
$\theta$, so that $L_{\rm tail}$ stays uniformly positive along the entire
segment from $x_0$ to $x_*$. Writing $f_\theta(x)=\det(\lambda I-M_x)$ with
$\lambda=e^{2\pi\ii\theta}$, we first compute $f_\theta'(x)/f_\theta(x)$ near
$x_*$ in the normalization of the relative Hill formula \cite{HOW}. Let
$\Phi(t)$ be the normalized fundamental matrix of \eqref{eq:whole-shape} in
position and velocity coordinates, and set $P_0=(I,0)$ and
$B_0=\binom{0}{I}$. Variation of constants gives
\[
 M'_x=M_x\int_0^{2\pi}\Phi(t)^{-1}B_0\q(t)^{-1}P_0\Phi(t)\dd t.
\]
Because the inhomogeneous forcing enters the first-order state equation as
$-B_0$ times the right-hand side and $P_0B_0=0$, the Green kernel of the
physical differential operator in \eqref{eq:whole-shape} relative to Lebesgue
measure is continuous on the diagonal, where it equals
$G_x(t,t)=-P_0\Phi(t)(\lambda I-M_x)^{-1}M_x\Phi(t)^{-1}B_0$. Jacobi's formula
therefore yields the relative trace identity
\begin{equation}\label{eq:relative-hill-trace}
 \frac{f_\theta'(x)}{f_\theta(x)}
 =\int_0^{2\pi}\q^{-1}\tr G_x(t,t)\dd t
 =\operatorname{Tr} T(\theta,x)^{-1}.
\end{equation}
Indeed, writing $L_x=-\partial_t^2+2\ii Z\partial_t+(H+xI)/\q$, the map
\[
 (Uw)(\phi)=\sqrt{\frac{2\pi}{\eta}}\,
            e^{-\ii\theta\phi} w(t(\phi)),\qquad
 \frac{\dd\phi}{2\pi}=\frac{\eta}{2\pi}\,\q(t)^{-1}\dd t,
\]
defines a unitary isomorphism from the twisted space $L^2(\q^{-1}\dd t)$ onto
periodic $L^2(\dd\phi/(2\pi))$ (satisfying
$\|Uw\|^2=\int_0^{2\pi}|w(t)|^2\q(t)^{-1}\dd t$) that conjugates $\q L_x$ to
$T(\theta,x)$, while $(\q L_x)^{-1}$ has integral kernel $G_x$ relative to
$\q^{-1}\dd t$. Coercive comparison with $1-\partial_\phi^2$ shows that
$T(\theta,x)^{-1}$ is trace class, and inserting Fourier Fej\'er means in the
periodic representation shows that its trace equals the diagonal integral of
its continuous kernel.

Wherever $T(\theta,x)$ is invertible, the square completion gives
$T^{-1}=XS^{-1}X^*+Q L_{\rm tail}^{-1}Q$. Since $\partial_x S=X^*X$,
integrating \eqref{eq:relative-hill-trace} from $x_*$ to $x$ yields
\begin{equation}\label{eq:relative-graph}
 \frac{f_\theta(x)}{f_\theta(x_*)}=
 \frac{\det S(\theta,x)}{\det S(\theta,x_*)}
 \det\nolimits_F\big(I+(x-x_*)L_{\rm tail}(\theta,x_*)^{-1}\big),
\end{equation}
where the Fredholm determinant is the convergent product over the eigenvalues
of the positive trace-class operator $L_{\rm tail}(\theta,x_*)^{-1}$. Identity
\eqref{eq:relative-graph} holds initially near $x_*$ and extends across zeros
of $\det S$ by analyticity. Near $(\theta_0,x_0)$, the tail determinant is
nonzero because $L_{\rm tail}$ is invertible, and it is real-analytic in
$\theta$: conjugating the tail operator by $(1-\partial_\phi^2)^{-1/2}$
produces a bounded coercive operator with analytic inverse flanked by two
Hilbert--Schmidt factors, so $L_{\rm tail}^{-1}$ depends analytically on
$\theta$ in trace norm. This proves \eqref{eq:same-characteristic}. The order
in $\theta$ follows because $\theta\mapsto e^{2\pi\ii\theta}$ has nonvanishing
derivative; at fixed $\theta$, the strict positivity $\partial_x S>0$ makes
the compression to $\ker S$ nondegenerate, so taking a Schur complement of
the invertible block on $(\ker S)^\perp$ gives vanishing order exactly
$\dim\ker S$.
\end{proof}

\begin{remark}[Lagrangian graphs]\label{rem:fourier-graphs}
On $V\oplus V$ with $V=PL^2$ and symplectic form
$\omega((u,p),(v,h))=u^*h-p^*v$, the graphs of $A_{\rm c}$ and $B_{\rm c}$
together with the vertical subspace $0\oplus V$ are Lagrangian, and the
vertical displacement from the graph of $B_{\rm c}$ to that of $A_{\rm c}$
recovers $u^*Sv$ via $\omega\bigl((u,B_{\rm c}u),(0,(A_{\rm c}-B_{\rm c})v)\bigr)=u^*Sv$.
Remark~\ref{rem:triple} gives the parallel real graph interpretation for the
complex-time matching pencil.
\end{remark}

\section{The periodic determinant and the real instability criterion}\label{sec:periodic-general}

\subsection{Periodic Hardy splitting and Jordan structure}
At $\theta=0$, the neighboring couplings $-ek(k+1)I/2$ in
\eqref{eq:physical-jacobi} vanish on both sides of the constant Fourier mode.
Writing $T_+(x)$ and $T_-(x)$ for the restrictions of $T(0,x)$ to the strictly
positive and strictly negative Fourier modes in $\phi$, we obtain the exact
orthogonal decomposition (valid even when $H+xI$ is singular)
\begin{equation}\label{eq:hardy-exact}
 T(0,x)=(H+xI)\oplus T_+(x)\oplus T_-(x).
\end{equation}
In particular, for a real physical system where $Z=\ii J$ ($J^T=-J$, $J^2=-I$)
and $H=R=R^T$, complex conjugation interchanges the two half-lines, so the real
periodic Morse index satisfies
\begin{equation}\label{eq:morse-parity}
 n_-(L_{R,e})=n_-(R)+2n_-(T_+(0)),
\end{equation}
where $n_-(T_+(0))$ is the complex negative index on the positive half-line
(the factor two accounting for the two conjugate halves of the complexified
form, whose total negative index equals the real index).

\begin{proposition}[Factorization of the periodic characteristic value]\label{prop:det-split}
Put $f(x)=\det(I-M_x)$ and choose $x_*\in\RR$ large enough that
$T(0,x_*)>0$. For each sign, let $\{\mu_j^\pm\}$ list the eigenvalues of
$T_\pm(x_*)$ with multiplicity and set
\[
 d_\pm(x)=\prod_j\left(1+\frac{x-x_*}{\mu_j^\pm}\right).
\]
Then for all $x\in\CC$,
\begin{equation}\label{eq:det-split}
 \frac{f(x)}{f(x_*)}
  =\frac{\det(H+xI)}{\det(H+x_*I)}d_+(x)d_-(x).
\end{equation}
For a real system one has $d_+=d_-$, whereas for a complex block the two
factors need not coincide. At any real $x$, the order of a zero of $f$ as a
function of the \emph{stiffness shift} equals the geometric multiplicity of the
periodic multiplier.
\end{proposition}

\begin{proof}
At $\theta=0$, the orthogonal decomposition \eqref{eq:hardy-exact} splits the
relative trace identity \eqref{eq:relative-hill-trace} into the logarithmic
derivatives of the three factors on the right-hand side of
\eqref{eq:det-split}. Both infinite products converge locally uniformly in
$x\in\CC$ because $\sum_j 1/\mu_j^\pm<\infty$ by elliptic comparison on each
half-line; since both sides of \eqref{eq:det-split} are entire and equal to
$1$ at $x_*$, they coincide everywhere. At any real zero, each vanishing
spectral factor contributes its kernel dimension, and in real coordinates
conjugation identifies the two half-line spectra at $x_*$, forcing $d_+=d_-$.
\end{proof}

\begin{proposition}[Signed periodic current and Jordan structure at $+1$]\label{thm:current-split}
Suppose $H$ is invertible and put $E_1=\{v:Hv=v,\ HZv=Zv\}$. Then the return
$M_0$ has exactly $\dim_{\CC}E_1$ Jordan blocks of size two at $+1$, while all
remaining blocks at $+1$ have size one. On the reducing complement of $E_1$,
the inner and outer periodic solution spaces are current-orthogonal and carry
positive and negative definite conserved currents, respectively.
\end{proposition}
\begin{proof}
Because $H$ is invertible, the constant block in \eqref{eq:hardy-exact} has
trivial kernel, while the kernels of $T_+(0)$ and $T_-(0)$ consist of periodic
solutions extending holomorphically to the interior and exterior of the unit
circle, respectively. Indeed, analytic ODE theory extends any periodic solution
to a neighborhood of $|z|=1$, and one-sided Fourier support in
$\zeta=(z-r)/(1-rz)$ gives holomorphy across the corresponding closed disk or
exterior; moreover, the absence of a constant term in $\zeta$ forces the inner
solution to vanish at $z=r$ and the outer solution to vanish at $z=r^{-1}$ (or
at $0$ and $\infty$ when $r=0$).

With the orthogonal projections $P_\pm=(I\pm Z)/2$, define in physical time
\[
 F_w=\q^{-1}e^{-\ii Zt}w
    =\frac{2(P_++z^2P_-)w}{ez^2+2z+e}.
\]
For $e>0$, the zero of $w$ at $r$ or $r^{-1}$ cancels the corresponding pole
of the denominator on each half; at $e=0$, the factors $z^{-1}P_++zP_-$
likewise preserve one-sided support on functions vanishing at $0$ or $\infty$.
Hence the physical Fourier series of $F_w$ in $z$ is supported in $k\ge 0$ on
the inner half and in $k\le 0$ on the outer half. Substituting $w=\q e^{\ii Zt}F_w$
into the constant pairing \eqref{eq:common-current} and averaging over
$t\in[0,2\pi]$ yields
\begin{equation}\label{eq:analytic-current}
 \tfrac12\left(\frac{1}{2\pi}\int_0^{2\pi}\q^{-2}\dd t\right)[u,w]_{\mathcal J}
       =\sum_{k\in\mathbb Z}k(F_u)_k^*(F_w)_k.
\end{equation}
The inner and outer spaces are therefore current-orthogonal, and equality to
zero in either definite sum of squares forces $F_w$ to be a constant vector
$v$. Using $\q''+\q=1$, the differential equation for $w=\q e^{\ii Zt}v$
reduces to $He^{\ii Zt}v=e^{\ii Zt}v$ for all $t$, which holds if and only if
$v\in E_1$. The orthogonal subspaces $E_1$ and $E_1^\perp$ reduce both $H$ and
$Z$; on $E_1^\perp$, the current pairing restricts nondegenerately to
$\ker(M_0-I)$, so $\operatorname{ran}(M_0-I)=\ker(M_0-I)^{\perp_{\mathcal J}}$
prevents any eigenvector there from admitting a Jordan predecessor. On $E_1$,
writing $w=e^{\ii Zt}h$ reduces the equation to the scalar system
$h''+h=h/\q$, whose linearly independent solutions $\q(t)$ and
$\q(t)\int_0^t\q(s)^{-2}\dd s$ generate one Jordan block of size two for each
basis direction in $E_1$, since the integral increases by a positive constant
over one period.
\end{proof}

\subsection{A signed determinant for the real return}
Write $p_*(z)=\det(zI-M_0)/(z-1)^{a_1}$, where $a_1$ is the algebraic
multiplicity of $+1$. The next proposition establishes the periodic parity law
for a real gyroscopic system on an even-dimensional position space.

\begin{proposition}[Parity after removing the periodic multiplier]\label{prop:stripped}
For every nonsingular real symmetric $R$, every real orthogonal complex
structure $J$, and every $0\le e<1$,
\begin{equation}\label{eq:physical-parity}
 \operatorname{sgn}p_*(1)=\operatorname{sgn}\det R
                =(-1)^{N_>(M_0)},
\end{equation}
where $N_>(M_0)$ counts the positive real multipliers greater than one with
their algebraic multiplicities.
\end{proposition}
\begin{proof}
Choose $x_*>0$ such that $R+x_*I>2I$. Since completing the square in the
kinetic term gives $\|w'+Jw\|^2-\|w\|^2$ and $\q\le 2$, the physical periodic
quadratic form remains strictly positive throughout a linear deformation to the
circular scalar endpoint $e=0$, $R+x_*I=cI$ with $c>2$. At that endpoint, the
periodic rotation $w=e^{-Jt}h$ reduces the equation to $h''=(c-1)h$, so each
real position direction contributes one negative factor
$(1-e^{2\pi\sqrt{c-1}})(1-e^{-2\pi\sqrt{c-1}})$ to $f(x_*)$; because the
position dimension is even, $f(x_*)>0$. Since $d_+=d_-$ in
\eqref{eq:det-split}, it follows that
$\operatorname{sgn}f(x)=\operatorname{sgn}\det(R+xI)$ at every regular real
$x$.

Set $E_{1,\RR}=\{v\in\RR^s:Rv=v,\ RJv=Jv\}$. Because $J^2=-I$, $J^T=-J$, and
$R^T=R$, the subspace $E_{1,\RR}$ is $J$-invariant, and both $E_{1,\RR}$ and
its orthogonal complement $W$ reduce $J$ and $R$; in particular, $W$ has even
real dimension and $R_\perp=R|_W$ satisfies $\det R_\perp=\det R$. Writing
$M_\perp(x)$ for the return on $W\oplus W$ and
$f_\perp(x)=\det(I-M_\perp(x))$, the same sign argument applied on $W$ gives
$\operatorname{sgn}f_\perp(x)=\operatorname{sgn}\det(R_\perp+xI)$ at every
regular real $x$. By Proposition~\ref{thm:current-split}, the unit-stiffness
block on $E_{1,\RR}$ contributes only the removed powers of $z-1$ at $x=0$, so
if $f_\perp(0)\neq 0$ then $f_\perp(0)=p_*(1)$ and the first equality in
\eqref{eq:physical-parity} holds immediately.

If $+1$ occurs in $M_\perp(0)$, it is semisimple by
Proposition~\ref{thm:current-split}. Let $K=\ker(M_\perp(0)-I)$ be its real
eigenspace. With respect to the preserved real symplectic form
$\Omega(u,w)=u^Tw'-u'^Tw+2u^TJw$, we have
$\operatorname{ran}(M_\perp(0)-I)=K^\Omega$, and semisimplicity implies
$K\cap K^\Omega=\{0\}$, so $K$ is symplectic and $k=\dim_{\RR}K$ is even.
Writing $M_\perp(0)=I_K\oplus N$ on $K\oplus K^\Omega$ and letting $Y$ be the
compression of $M_\perp'(0)$ to $K$ in this splitting, variation of constants
at fixed initial data gives
\[
 \Omega(v,Yu)=\int_0^{2\pi}\q^{-1}w_v(t)^T w_u(t)\dd t\qquad(u,v\in K).
\]
Since the right-hand side is positive definite and the nondegenerate real skew
matrix representing $\Omega|_K$ has positive determinant, we obtain $\det Y>0$.
Expanding $f_\perp(x)$ in the block splitting $K\oplus K^\Omega$ then yields
\[
 f_\perp(x)=(-x)^k\det Y\,\det(I-N)+O(x^{k+1}).
\]
Because $k$ is even, $R_\perp$ is nonsingular, and $p_*(1)=\det(I-N)$, passing
to small nonzero real $x$ gives
$\operatorname{sgn}p_*(1)=\operatorname{sgn}\det R_\perp=\operatorname{sgn}\det R$.
Finally, nonreal conjugate pairs and negative real multipliers contribute
positively to $p_*(1)$ through factors $|1-\lambda|^2$ and $1-\lambda>0$, respectively, while each real
multiplier greater than one contributes a negative factor, proving the second
equality in \eqref{eq:physical-parity}.
\end{proof}

\begin{corollary}[Newtonian planar shapes]\label{cor:newton-all}
For a collision-free planar central configuration of positive masses, let $R$
be the normalized constrained Hessian on the essential shape space. If $R$ is
invertible, then every essential periodic multiplier is semisimple for all
$0\le e<1$, with current signature $(d,d)$ on its complex eigenspace, where $2d=\nu_1(M_0)$;
moreover, if $n_-(R)$ is odd, then $M_0$ admits a positive real multiplier
greater than one at every $0\le e<1$, including periodic resonances.
\end{corollary}
\begin{proof}
Work in the mass inner product with center of mass at the origin, normalize
the configuration $a_0$ by $\nabla_m U(a_0)=-\lambda_c a_0$ for the positive
Newtonian potential $U$, and set $L=D\nabla_m U(a_0)$. Translation invariance,
homogeneity, and rotational equivariance give $L=0$ on translations,
$La_0=2\lambda_c a_0$, and $LJa_0=-\lambda_c Ja_0$, so their mass-orthogonal
complement $W$ reduces both $L$ and $J$, with normalized constrained Hessian
$R=I+L|_W/\lambda_c$. Writing a variation as $s e^{Jt}w(t)$ in terms of the
Kepler radius factor $s=h^2/(\lambda_c\q)$ and true anomaly $t$ (where
$h=s^2\dot t$ is constant and $\ddot s-s\dot t^2=-\lambda_c/s^2$) turns the
rotated physical Hessian $s^{-3}L$ into the equation $w''+2Jw'=\q^{-1}Rw$ on
$W$, realizing the essential planar return of the Meyer--Schmidt decomposition
\cite{MS}.

For any real $v\in W$, summing the pairwise edge Hessians in the orthogonal
directions $v$ and $Jv$ yields
\[
 \langle Lv,v\rangle_m+\langle LJv,Jv\rangle_m
       =\sum_{i<j}\frac{m_i m_j}{|a_{0,i}-a_{0,j}|^3}|v_i-v_j|^2,
\]
since each edge matrix has the form $3nn^T-I$ and the squared projections of
$v_i-v_j$ and its quarter turn onto the unit edge vector $n$ sum to
$|v_i-v_j|^2$. Because $v\in W$ is orthogonal to translations, this sum is
strictly positive whenever $v\neq 0$, so the $J$-invariant subspace in
$\ker(R-I)$ vanishes. Proposition~\ref{thm:current-split} therefore implies
semisimplicity at $+1$, complex conjugation interchanges the two definite
halves to give signature $(d,d)$, and \eqref{eq:physical-parity} produces a
positive real multiplier greater than one whenever $n_-(R)$ is odd. Together
with Proposition~\ref{prop:stripped}, this proves
\hyperref[main:parity]{Theorem~A}.
\end{proof}

\subsection{The characteristic polynomial from a reversible half-period}
\label{subsec:reversible-half}

For a real shape equation with an orthogonal involution $Q$ satisfying
$QR=RQ$ and $QJ=-JQ$, time reflection acts on position and velocity by
$K=\diag(Q,-Q)$. This is an antisymplectic involution for the form
$\Omega$ above. Since $q_e(2\pi-t)=q_e(t)$, reflection gives
$M_0=KT^{-1}KT$, where $T$ is the transition from $0$ to $\pi$.
The half-period factorization is standard \cite[Section~3.1]{PWreversible}.
We express the reduced characteristic polynomial of a reversible return
\cite{FKM} directly in terms of the blocks of $T$.

\begin{proposition}[Reversible half-period polynomial and ranks]
\label{prop:reversible-half}
Let $T$ be a symplectic map on a real symplectic space of dimension $2d$,
let $K$ be an antisymplectic involution, and set $P=KT^{-1}KT$.
Write $V_\pm=\ker(K\mp I)$ and, relative to $V_+\oplus V_-$, put
\[
 T=\begin{pmatrix}A&B\\ L&D\end{pmatrix},\qquad
 G(z)=\det\begin{pmatrix}A&B\\ zL&D\end{pmatrix}
     =\sum_{k=0}^{d}g_k z^k,\qquad
 \widetilde G(s,t)=\sum_{k=0}^{d}g_k s^{d-k}t^k.
\]
Then $G(1)=1$ and
\begin{equation}\label{eq:reversible-half-char}
 \det(\lambda I-P)=\widetilde G\bigl((1-\lambda)^2,(1+\lambda)^2\bigr).
\end{equation}
The return is spectrally stable exactly when all projective roots of
$\widetilde G$, in the coordinate $z=t/s$, belong to
$(-\infty,0]\cup\{\infty\}$. A finite nonzero root of multiplicity $m$
gives two reciprocal multipliers, each of algebraic multiplicity $m$
and geometric multiplicity $\dim_{\CC}\ker\left(\begin{smallmatrix}A&B\\
zL&D\end{smallmatrix}\right)$. At the endpoints one has
\[
\begin{array}{c|cc}
 &\text{algebraic multiplicity}&\text{eigenspace}\\ \hline
 -1&2\operatorname{ord}_{0}G&\ker A\oplus\ker D\\
 +1&2(d-\deg G)&\ker L\oplus\ker B .
\end{array}
\]
The endpoint kernels in the table use the original splitting
$V_+\oplus V_-$. The return has bounded powers exactly when it is
spectrally stable and each geometric multiplicity equals its algebraic
multiplicity.
\end{proposition}

\begin{proof}
The eigenspaces of $K$ are complementary Lagrangian subspaces, each of
dimension $d$. Put $a=1-\lambda$, $b=1+\lambda$. The identity
\[
 KTK(P-\lambda I)=T-\lambda KTK
   =\begin{pmatrix}a A&b B\\ b L&a D\end{pmatrix}
\]
preserves kernels and determinants since $\det(KTK)=1$ and $2d$ is even.
For $ab\ne0$, put $S=\diag(I,(b/a)I)$. Conjugating the displayed
matrix divided by $a$, with $S$ on the left and $S^{-1}$ on the right,
gives $\left(\begin{smallmatrix}A&B\\ (b/a)^2L&D\end{smallmatrix}\right)$.
Thus its determinant before division is $a^{2d}G(b^2/a^2)$.
As the homogenized expression is polynomial in $\lambda$, the identity
holds at $\lambda=\pm1$ as well. The block identity at these two values
gives the endpoint eigenspaces.

The unit circle is exactly the inverse image of
$(-\infty,0]\cup\{\infty\}$ under $z=((1+\lambda)/(1-\lambda))^2$,
since $(1+\lambda)/(1-\lambda)$ is purely imaginary exactly when
$|\lambda|=1$ for $\lambda\ne 1$.
The map has local degree two at $\lambda=\pm1$ and is unramified elsewhere. The value $z=1$ is excluded
because $G(1)=1$. These facts give the algebraic multiplicities, including
the multiplicity $d-\deg G$ of the projective root at infinity. The same
invertible block transformations give the geometric multiplicities away
from the endpoints. The corresponding chain equations are
$(T-\lambda_0 KTK)v_j=KTK v_{j-1}$ if individual Jordan block sizes are
needed. By symplectic reciprocity, a multiplier off the unit circle implies
one of modulus greater than one. On the unit circle, the Jordan form gives
the boundedness criterion.
\end{proof}

For the comparison with the reduced return, let $P_{++}$ be the diagonal
block of $P$ on $V_+$. With $X=(T^{-1})_{+-}L$, block multiplication gives
$P_{++}=I-2X$, and left multiplication of the matrix defining $G$ by $T^{-1}$
gives a block triangular matrix with diagonal blocks $I+(z-1)X$ and $I$.
Consequently,
\[
 \widetilde G(s,t)=2^{-d}\det\bigl((s+t)I+(s-t)P_{++}\bigr).
\]
In particular, if $c_1,\ldots,c_d$ are the roots of
$\det(cI-P_{++})$, the full characteristic polynomial is
$\prod_{j=1}^{d}(\lambda^2-2c_j\lambda+1)$. The identity $P^{-1}=KPK$
and symplecticity also give $P_{--}=P_{++}^{\,T}$ in dual Lagrangian bases.

For a reversible shape equation, write $n=\deg G$. After removing the periodic
factor as in Proposition~\ref{prop:stripped}, \eqref{eq:reversible-half-char}
gives $p_*(1)=4^n g_n$. Thus the odd-index instability criterion is read from the
sign of the highest nonzero coefficient of $G$, even at a periodic resonance.
Positive real multipliers greater than one correspond to roots $z>1$.
At the other endpoint, \eqref{eq:reversible-half-char} gives
\[
 \det(I+P)=2^{2d}\det A\,\det D .
\]
An antiperiodic solution starting in $V_+$ must end in $V_-$ at $\pi$, and a
solution starting in $V_-$ must end in $V_+$. The two kernels $\ker A$ and
$\ker D$ are therefore the two mixed reflection-sector kernels. If the root at $z=0$ is simple, then $c=-1$ is simple in each
diagonal block of $P$. Indeed, the diagonal blocks of $T^{-1}T=I$ give
$P_{++}+I=2(T^{-1})_{++}A$ and $P_{--}+I=2(T^{-1})_{--}D$.
The sector kernels thus lie in the respective block eigenspaces, so each
singular sector has nullity one. When both sectors are singular,
$-1$ is semisimple, while a zero in just one sector gives a single
Jordan block of size two.

The same representation gives a one-point test for spectral instability.
Indeed, spectral stability forces the eigenvalues $x_j$ of $X$ into $[0,1]$,
since $c_j=1-2x_j$. Hence
$G(\ii)=\prod_j(1-x_j+\ii x_j)$ satisfies
$2^{-d/2}\le |G(\ii)|\le1$. Either strict violation detects a multiplier off
the unit circle.

\section{Algebraic matching for two local monodromies}\label{sec:theory}

Throughout this section we work in the complex fiber $V$ at the common base
point $z_b$. A simultaneous change of basis conjugates the pencil
\eqref{eq:master} below and carries every kernel, image, and quotient space to
its counterpart. For symplectic isometries, de~Seguins~Pazzis \cite{DSP}
studies which isometries are products of two unipotents of index at most two;
in our equations the two local singularities prescribe the factors, and
analytic continuation along $(r,0)$ determines their relative position.

\begin{theorem}[Complex Cayley pencil and reciprocity criteria]\label{thm:cayley}
Let $V$ be a finite-dimensional complex vector space, let
$U_0,U_r\in\mathrm{GL}(V)$ have no eigenvalue $-1$, and put $M=U_0U_r$. Define
the Cayley generators
\[
 A=-2\ii(U_0-I)(U_0+I)^{-1},\qquad B=-2\ii(U_r-I)(U_r+I)^{-1}.
\]
\begin{enumerate}
\item For every $\lambda\in\CC$, multiplication by $I-\ii B/2$ identifies
$\ker\mathcal C(\lambda)$ with $\ker(M-\lambda I)$, where
\begin{equation}\label{eq:master}
 \mathcal C(\lambda)=(1-\lambda)(I-AB/4)+\tfrac{\ii}{2}(1+\lambda)(A+B).
\end{equation}
In particular, the antiperiodic equation reduces to $ABx=4x$. At $\lambda=1$,
put $E=\ker A\cap\ker B$ and $W_\cap=\im A\cap\im B$, and define the
obstruction map
\[
 \mathcal O:W_\cap\longrightarrow V/(\ker A+\ker B),\qquad
 \mathcal O x=[y_A+y_B]\quad(Ay_A=By_B=x).
\]
Then $\nu_1(M)=\dim E+\dim\ker\mathcal O$; moreover, when $A^2=B^2=0$, one has
$W_\cap\subset E$ and $\dim E\le\nu_1(M)\le 2\dim E$.

\item Suppose in addition that $\dim_{\CC}V=4$, $A^3=0$, and $B^2=0$, and set
\[
 \Delta(M)=\tr M-\tr M^{-1},\qquad \chi_M(\lambda)=\det(\lambda I-M).
\]
Then $\Delta(M)=-\ii\tr(A^2B)$, and the following assertions hold.
\begin{enumerate}
\item If $\lambda_0\neq 1$ satisfies $\nu_{\lambda_0}(M)\ge 2$, then with
$\rho_0=2-\lambda_0-\lambda_0^{-1}$,
\[
 \chi_M(\lambda)=\bigl(\lambda^2-(2-\rho_0)\lambda+1\bigr)^2.
\]
Consequently, $\Delta(M)\neq 0$ implies $\nu_\lambda(M)\le 1$ for every
$\lambda\neq 1$. Furthermore, at any root $\lambda=\pm 1$, one has
$\chi_M'(\lambda)=-\Delta(M)$.
\item If $A^2=0$ and $M=M(0)$ lies on a $C^1$ curve of matrices with constant
determinant such that $\frac{\mathrm d}{\mathrm d\varepsilon}\Delta(M(\varepsilon))\big|_0\neq 0$,
then $\nu_\lambda(M)\le 1$ for every $\lambda\neq 1$.
\end{enumerate}
\end{enumerate}
\end{theorem}

\begin{proof}[Proof of Theorem~\ref{thm:cayley}\textup{(1)}]
By definition of the Cayley generators $A$ and $B$, the matrices
$I-\ii A/2$ and $I-\ii B/2$ are invertible and satisfy
\[
 U_0=(I+\ii A/2)(I-\ii A/2)^{-1},\qquad
 U_r=(I+\ii B/2)(I-\ii B/2)^{-1}.
\]
Clearing these invertible denominators on the left and right of $M-\lambda I$
yields
\[
 (I-\ii A/2)(M-\lambda I)(I-\ii B/2)
 =(I+\ii A/2)(I+\ii B/2)-\lambda(I-\ii A/2)(I-\ii B/2)
 =\mathcal C(\lambda).
\]
Hence multiplication by $I-\ii B/2$ maps $\ker\mathcal C(\lambda)$
isomorphically onto $\ker(M-\lambda I)$, while $\det\mathcal C(\lambda)$
differs from $\det(M-\lambda I)$ by a nonzero factor independent of
$\lambda$. For the differential systems of \eqref{eq:whole-shape}, whenever
the complex-time factors are defined as in Section~\ref{subsec:setup},
Proposition~\ref{prop:graph-det} therefore identifies the algebraic
multiplicity read from $\mathcal C(\lambda)$ with the phase order of the
finite Schur determinant, including at $\lambda=\pm 1$. Setting $\lambda=-1$
cancels the sum $A+B$ and reduces the kernel condition to $ABx=4x$, whereas
setting $\lambda=1$ leaves $\ii(A+B)$.

At $\lambda=1$, let $x\in W_\cap=\im A\cap\im B$ and choose $y_A,y_B\in V$ with
$Ay_A=By_B=x$. Any alternative choices alter $y_A+y_B$ by an element of
$\ker A+\ker B$, so $\mathcal O:W_\cap\to V/(\ker A+\ker B)$ is well-defined
and linear. Moreover, $\mathcal O x=0$ holds precisely when
$y_A+y_B=k_A+k_B$ for some $k_A\in\ker A$ and $k_B\in\ker B$, in which case
$v=y_A-k_A=-y_B+k_B$ satisfies $Av=x$ and $Bv=-x$; conversely, any
$v\in\ker(A+B)$ with $Av=x$ supplies the preimages $y_A=v$ and $y_B=-v$.
Therefore the sequence
\begin{equation}\label{eq:seq}
 0\longrightarrow E\longrightarrow\ker(A+B)\xrightarrow{\ A\ }
 \ker\mathcal O\longrightarrow 0
\end{equation}
is exact, giving $\nu_1(M)=\dim E+\dim\ker\mathcal O$. Finally, when
$A^2=B^2=0$, one has $\im A\subset\ker A$ and $\im B\subset\ker B$, so
$W_\cap\subset E$ and $\dim E\le\nu_1(M)\le 2\dim E$. (If in addition
$\rank A=\rank B=\frac12\dim V$, then $W_\cap=E$ and
$\nu_1(M)=2\dim E-\rank\mathcal O$.)
\end{proof}

To prove part~(2) of Theorem~\ref{thm:cayley}, we first compress the nonperiodic eigenvalue
problem onto the at most two-dimensional image of the pole generator $B$.

\begin{lemma}[Compression onto the image of the pole generator]\label{lem:ranktwo}
Under the four-dimensional hypotheses of Theorem~\ref{thm:cayley}\textup{(2)}, one has
$\nu_\lambda(M)\le\rank B\le 2$ for every $\lambda\neq 1$. If $\rank B=2$,
choose a full-rank factorization $B=E_BF_B^T$ over $\CC$ and set
$C_1=F_B^TAE_B$, $C_2=F_B^TA^2E_B$, and $d_\lambda=1-\lambda$. Then
$\ker(M-\lambda I)\simeq\ker\mathcal R_\lambda$ for every $\lambda\neq 1$,
where
\begin{equation}\label{eq:reduced}
 \mathcal R_\lambda=I+\frac{\lambda}{d_\lambda^2}C_1
 -\frac{\ii\lambda(\lambda+1)}{2d_\lambda^3}C_2,\qquad
 \det(M-\lambda I)=d_\lambda^4\det\mathcal R_\lambda.
\end{equation}
In particular, $\mathcal R_{-1}=I-C_1/4$ and $\tr(\adj(C_1)C_2)=0$.
\end{lemma}
\begin{proof}
Since $A^3=0$ and $B^2=0$, the Cayley formulas expand finitely as
$U_0=I+\ii A-A^2/2$ and $U_r=I+\ii B$. For $\lambda\neq 1$, the matrix
$U_0-\lambda I$ is invertible and $M-\lambda I=(U_0-\lambda I)+\ii U_0B$
differs from it by a perturbation of rank $\rank B$, which gives
$\nu_\lambda(M)\le\rank B\le 2$.

When $\rank B=2$, the identity $B^2=0$ forces $F_B^TE_B=0$. Setting
$y=F_B^Tv$ eliminates the four-dimensional state $v$ from $(M-\lambda I)v=0$
via the mutually inverse maps
\[
 v=-\ii(U_0-\lambda I)^{-1}U_0E_By,\qquad
 \bigl(I+\ii F_B^T(U_0-\lambda I)^{-1}U_0E_B\bigr)y=0.
\]
Using $(U_0-\lambda I)^{-1}U_0=d_\lambda^{-1}I-d_\lambda^{-2}\ii\lambda A
-d_\lambda^{-3}\frac{\lambda(\lambda+1)}{2}A^2$ together with $F_B^TE_B=0$
yields $\mathcal R_\lambda$ in \eqref{eq:reduced}, and the Sylvester
determinant identity $\det(I+XY)=\det(I+YX)$ gives
$\det(M-\lambda I)=d_\lambda^4\det\mathcal R_\lambda$. Because a four-by-four
nilpotent matrix with $A^3=0$ has at most one Jordan chain of length three, we
have $\rank C_2\le 1$ and hence $\det C_2=0$. Expanding
$d_\lambda^4\det\mathcal R_\lambda$ therefore leaves only one potentially
singular term at $\lambda=1$, namely
\[
 -\frac{\ii\lambda^2(\lambda+1)}{2d_\lambda}\tr(\adj(C_1)C_2).
\]
Since $\det(M-\lambda I)$ is a polynomial in $\lambda$, this residue vanishes.
\end{proof}

\begin{proof}[Proof of Theorem~\ref{thm:cayley}\textup{(2)}]
Expanding $M=(I+\ii A-A^2/2)(I+\ii B)$ and
$M^{-1}=(I-\ii B)(I-\ii A-A^2/2)$ yields
\[
 \tr M=4-\tr(AB)-\tfrac{\ii}{2}\tr(A^2B),\qquad
 \tr M^{-1}=4-\tr(AB)+\tfrac{\ii}{2}\tr(A^2B),
\]
whence $\Delta(M)=-\ii\tr(A^2B)$. Because both factors are unipotent,
$\det M=1$. Writing $t_1=\tr M$ and $\tilde t_1=\tr M^{-1}$, the characteristic
polynomial takes the form
$\chi_M(\lambda)=\lambda^4-t_1\lambda^3+\chi_2\lambda^2-\tilde t_1\lambda+1$,
so evaluating $\chi_M'(\pm 1)$ at a root $\chi_M(\pm 1)=0$ eliminates $\chi_2$
and gives $\chi_M'(\pm 1)=-(t_1-\tilde t_1)=-\Delta(M)$.

Suppose now that $\lambda_0\neq 1$ satisfies $\nu_{\lambda_0}(M)\ge 2$. By
Lemma~\ref{lem:ranktwo}, we must have $\rank B=2$ and $\mathcal R_{\lambda_0}=0$,
so \eqref{eq:reduced} expresses $C_1=\rho_0 I+\alpha_0 C_2$ for some scalar
$\alpha_0\in\CC$, where $\rho_0=-d_{\lambda_0}^2/\lambda_0=2-\lambda_0-\lambda_0^{-1}\neq 0$.
In dimension two the adjugate is linear, and $\rank C_2\le 1$ implies
$\adj(C_2)C_2=0$; substituting $C_1$ into $\tr(\adj(C_1)C_2)=0$ therefore
yields $\rho_0\tr C_2=0$, so $\tr C_2=\det C_2=0$ and $C_2$ is nilpotent. For
every $\lambda\neq 1$, the matrix $\mathcal R_\lambda$ in \eqref{eq:reduced} is
consequently the scalar $(1+\lambda\rho_0/d_\lambda^2)I$ plus a multiple of the
nilpotent matrix $C_2$. Taking its determinant and extending by continuity to
$\lambda=1$ gives
\[
 \chi_M(\lambda)=d_\lambda^4\left(1+\frac{\lambda\rho_0}{d_\lambda^2}\right)^2
 =\bigl(\lambda^2-(2-\rho_0)\lambda+1\bigr)^2.
\]
This polynomial is reciprocal, forcing $\Delta(M)=0$ and establishing
part~(2a) (in particular, $\lambda_0=-1$ forces $\chi_M(\lambda)=(\lambda+1)^4$).

For part~(2b), assume $A^2=0$ and suppose for contradiction that some
$\lambda_0\neq 1$ has $\nu_{\lambda_0}(M)\ge 2$. Then $\rank B=2$, $C_2=0$,
and $\mathcal R_{\lambda_0}=0$ yields $C_1=\rho_0 I$, so
$BAE_B=E_BF_B^TAE_B=\rho_0 E_B$. Together with $BE_B=0$ and $\rho_0\neq 0$,
this implies that the columns of $(E_B,AE_B)$ form a basis of $V$, in which
\[
 A=\begin{pmatrix}0&0\\I&0\end{pmatrix},\qquad
 B=\begin{pmatrix}0&\rho_0 I\\0&0\end{pmatrix},\qquad
 M+M^{-1}=(2-\rho_0)I.
\]
Along any $C^1$ curve $M(\varepsilon)$ passing through $M=M(0)$, differentiating
$\tr(M-M^{-1})$ and using cyclicity of the trace gives
\begin{equation}\label{eq:tangent}
 \dot\Delta=\tr\bigl((M+M^{-1})M^{-1}\dot M\bigr).
\end{equation}
Since $\det M(\varepsilon)$ is constant, Jacobi's formula gives
$\tr(M^{-1}\dot M)=0$, and substituting $M+M^{-1}=(2-\rho_0)I$ into
\eqref{eq:tangent} forces $\dot\Delta=0$ at $\varepsilon=0$, a contradiction.
\end{proof}

When $A^2=0$ and the common fixed space $E=\ker A\cap\ker B$ is
one-dimensional, a single left fixed covector also detects a two-dimensional
periodic eigenspace even without reflection symmetry.

\begin{corollary}[Periodic multiplicity two from a fixed covector]\label{cor:fixed-covector}
Under the hypotheses of Theorem~\ref{thm:cayley}\textup{(2)}, suppose $A^2=0$, $\dim E=1$,
and there exist $g\in E$ and a complex linear functional $f\in V^*$ such that
$fM=f$ and $f(g)\neq 0$. Then $\nu_1(M)=2$. In particular, if $M$ preserves a
nondegenerate Hermitian form $[\,\cdot,\cdot\,]$ (conjugate-linear in its first
argument), any $g\in E$ with $[g,g]\neq 0$ yields $\nu_1(M)=2$.
\end{corollary}
\begin{proof}
Since $A^2=0$, we have $\Delta(M)=0$; moreover, every nonzero $g\in E$ is fixed
by both $U_0$ and $U_r$, so $+1$ is an eigenvalue of $M$ with
$\chi_M'(1)=-\Delta(M)=0$ and hence algebraic multiplicity at least two. By
Theorem~\ref{thm:cayley}(1), $\nu_1(M)\le 2\dim E=2$. If $\nu_1(M)=1$, then the single Jordan
chain at $+1$ would place $g\in\im(M-I)$, contradicting
$f\bigl((M-I)v\bigr)=(fM-f)(v)=0$ for all $v\in V$. In the Hermitian case,
$f(v)=[g,v]$ satisfies $f(Mv)=[g,Mv]=[M^{-1}g,v]=[g,v]=f(v)$.
\end{proof}

\begin{remark}[Realization of a double antiperiodic root]\label{rem:double-ap}
The trace and variational tests in Theorem~\ref{thm:cayley}(2) genuinely probe the relative
position of the two local factors. On $\CC^2\oplus\CC^2$, the pair
\[
 A=\begin{pmatrix}0&0\\I&0\end{pmatrix},\qquad
 B=\begin{pmatrix}0&4I\\0&0\end{pmatrix},\qquad
 M=(I+\ii A)(I+\ii B)=\begin{pmatrix}I&4\ii I\\\ii I&-3I\end{pmatrix}
\]
satisfies $A^2=B^2=0$ and $(M+I)^2=0$ with $\rank(M+I)=2$, realizing both
$\chi_M(\lambda)=(\lambda+1)^4$ and $\nu_{-1}(M)=2$.
\end{remark}

\begin{remark}[Symplectic and Hermitian forms of the periodic obstruction]\label{rem:forms}
When both local factors preserve a common nondegenerate pairing, the quotient
in \eqref{eq:seq} admits a concrete dual description. If $U_0$ and $U_r$
preserve a nondegenerate alternating bilinear form $\omega$, then $A$ and $B$
are Hamiltonian, so $\ker A=(\im A)^\omega$ and $\ker B=(\im B)^\omega$, and
$\omega$ identifies $V/(\ker A+\ker B)$ with $W_\cap^*$. Defining symmetric
forms $g_A(Ay,v)=\omega(y,v)$ on $\im A$ and $g_B(By,v)=\omega(y,v)$ on
$\im B$, the obstruction evaluates on $v\in W_\cap$ as
$\omega(\mathcal O x,v)=(g_A+g_B)(x,v)$, so $\ker\mathcal O$ is precisely the
radical of $(g_A+g_B)|_{W_\cap\times W_\cap}$.

Similarly, if both factors preserve a common nondegenerate Hermitian form
$[\,\cdot,\cdot\,]$, then $A$ and $B$ are self-adjoint for it (the factor
$-\ii$ converting the skew-adjoint Cayley transform into a self-adjoint one),
and the induced forms $g_A(Ay,v)=[y,v]$ and $g_B(By,v)=[y,v]$ are Hermitian
with the same radical characterization of $\ker\mathcal O$. In our polygonal
complex blocks below, Hermitian invariance holds for the \emph{product} $M$
along real time rather than for the individual based factors; the conserved
current of a common holomorphic solution then supplies both the covector in
Corollary~\ref{cor:fixed-covector} and, via time reflection, the signature
$(1,1)$ in Proposition~\ref{prop:positive-current}.
\end{remark}

\section{Real two-component systems and the Euler--Lagrange families}\label{sec:real}

We now apply the complex-time factorization to the real variational equation
\eqref{eq:real-op}. Although the physical period map lies in $\mathrm{Sp}(4,\RR)$,
its complex-time local factors $U_0=I+\ii A$ and $U_r=I+\ii B$ are complex
symplectic and nonreal whenever their generators are nonzero in the real fiber
used below. The Hamiltonian symmetry at the origin forces $A^2=0$, which is
the exact algebraic distinction between the real system \eqref{eq:real-op} and
the complex block \eqref{eq:complex-op}.

\subsection{Local square-zero generators and the symmetric pencil}\label{subsec:real-pencil}

\begin{lemma}[Local generators and the highest-exponent solution]\label{lem:real-local}
For \eqref{eq:real-op} with $0<e<1$, there exists a common real symplectic
fiber at $z_b\in(r,0)$ in which $U_0=I+\ii A$ and $U_r=I+\ii B$, where $A,B$
are real Hamiltonian endomorphisms satisfying $A^2=B^2=0$ and
$\rank B=\rank D$. Moreover, if $\rank A=2$, then $\im A$ contains a real
vector $\phi=Ay$ with $\omega(y,\phi)=0$ whose position component, in the
coordinates below, has the expansion $x_\phi(z)=z^2(1,1)^T+O(z^3)$.
\end{lemma}
\begin{proof}
An $SO(2)$ rotation preserves $J$ and diagonalizes
$D=\diag(\alpha_1,\alpha_2)$ with $\alpha_1\le\alpha_2$. Setting
$u=\diag(1,\ii)x$ and $X=\bigl(\begin{smallmatrix}0&1\\1&0\end{smallmatrix}\bigr)$
transforms \eqref{eq:real-op} into $\vang^2 x-2X\vang x+\q^{-1}Dx=0$, which in
the state variable $W=(Z(\vang x-Xx),x)^T$ reads
\begin{equation}\label{eq:fuchs-real}
 \vang W=\begin{pmatrix}-X&Z(I-D/\q)\\Z&X\end{pmatrix}W.
\end{equation}
The coefficient matrix is real along $(r,0)$ and Hamiltonian with respect to
the canonical symplectic form $\omega(v,y)=v^TJ_4y$, where
$J_4=\bigl(\begin{smallmatrix}0&-I\\I&0\end{smallmatrix}\bigr)$. In ordinary
$z$-derivative form, the residues at $0$ and $r$ are
\begin{equation}\label{eq:res-real}
 R_0=\begin{pmatrix}-X&Z\\Z&X\end{pmatrix},\qquad
 R_r=\begin{pmatrix}0&-ZD/\eta\\0&0\end{pmatrix},
\end{equation}
where we used $r\q'(r)=\eta$.

At $r$, all eigenvalues of $R_r$ vanish, so a real symplectic analytic gauge
$G_r(z-r)$ with $G_r(0)=I$ reduces the local system to $R_r$. At $0$, the
residue satisfies $R_0^3=4R_0$, so its exponent spaces $V_{-2},V_0,V_2$ have
dimensions $1,2,1$. Successive real symplectic changes $\exp(z^jG_j)$ with
$G_j\in\mathfrak{sp}(4,\RR)$ solve the homological equation on the nonresonant
subspace of $\mathfrak{sp}(4,\RR)$ at each order, leaving a resonant normal
form $R_0+z^2N_2+z^4N_4$ with $N_2,N_4\in\mathfrak{sp}(4,\RR)$ and
$[R_0,N_j]=jN_j$. By Lemma~\ref{lem:gauge-rec} in Appendix~\ref{app:gauge},
this formal symplectic gauge converges and retains $G^TJ_4G=J_4$.

Put $N=N_2+N_4$. Because $(-z)^{R_0}N_j(-z)^{-R_0}=z^jN_j$ for even $j$, the
local fundamental matrix at the origin may be taken as
\[
 G(z)(-z)^{R_0}\exp\bigl(N\log(-z)\bigr),
\]
which is real and symplectic on $(r,0)$. Near $r$, we take
$G_r(z-r)\exp\bigl(R_r\log(z-r)\bigr)$. Continuing both fundamental matrices to
$z_b$ conjugates $2\pi N$ and $2\pi R_r$ to real Hamiltonian endomorphisms
$A,B$ in the common symplectic fiber at $z_b$, with $\rank B=\rank R_r=\rank D$.

To show $N^2=0$, fix nonzero generators $e_\pm\in V_{\pm 2}$. Since $N$ raises
exponents by $2$ or $4$, there exist $v\in V_0$, a scalar $\alpha_0$, and a
linear functional $\ell_0$ on $V_0$ such that $Ne_-=v+\alpha_0 e_+$,
$Ny=\ell_0(y)e_+$ for $y\in V_0$, and $Ne_+=0$. Because $R_0$ is Hamiltonian,
$\omega$ pairs $V_{-2}$ nondegenerately with $V_2$ and makes $V_0$ orthogonal
to both. Applying the Hamiltonian identity $\omega(Ne_-,y)+\omega(e_-,Ny)=0$
then yields $\ell_0(y)=-\omega(v,y)/\omega(e_-,e_+)$, whence $\ell_0(v)=0$ and
consequently $N^2=0$. If $\rank N=2$, then $v\ne0$ since otherwise $\im N\subset V_2$. Nondegeneracy on $V_0$ now gives $\ell_0\ne0$;
choosing $y\in V_0$ with $\ell_0(y)=1$ gives $Ny=e_+$ and $\omega(y,e_+)=0$.
Transporting this pair to $z_b$ (and absorbing the factor $2\pi$ into the
preimage) yields $\phi=Ay$ with $\omega(y,\phi)=0$. By \eqref{eq:fuchs-real},
the position component of $e_+$ may be normalized to $(1,1)^T$, giving
$x_\phi(z)=z^2(1,1)^T+O(z^3)$.
\end{proof}

\begin{proposition}[Symmetric pencil and periodic intersection formula]\label{prop:shear}
Let $A,B$ be square-zero real Hamiltonian endomorphisms of a
$2d_V$-dimensional real symplectic space $(V_{\RR},\omega)$, and set
$M=(I+\ii A)(I+\ii B)$ on the complexification of $V_{\RR}$. On $L=\im A$,
define
\[
 g_A(Ay,x)=\omega(y,x),\qquad b_L(x,y)=\omega(x,By),\qquad C_L=AB|_L.
\]
Then $g_A$ is nondegenerate, both $g_A$ and $b_L$ are symmetric, and for every
$\lambda\neq 0,1$ one has (upon bilinear complexification)
\begin{equation}\label{eq:pencil}
 \ker(M-\lambda I)\simeq\ker(C_L-\rho I)=\radform(b_L+\rho g_A),\qquad
 \rho=2-\lambda-\lambda^{-1}.
\end{equation}
In particular, $\nu_\lambda(M)\le\min(\rank A,\rank B)$. Writing
$d_\lambda=1-\lambda$ and $r_A=\rank A$, the characteristic polynomial satisfies
\[
 \det(M-\lambda I)=d_\lambda^{\,2d_V-2r_A}
 \det\nolimits_L(d_\lambda^2 I+\lambda C_L).
\]
At $\lambda=1$, set $L_r=\im B$, $W_\cap=L\cap L_r$, and define $g_B$ on $L_r$
analogously to $g_A$. Then
\[
 \nu_1(M)=2d_V-\rank A-\rank B+\dim W_\cap
 +\dim\radform\bigl((g_A+g_B)|_{W_\cap\times W_\cap}\bigr).
\]
\end{proposition}
\begin{proof}
Since $\ker A=L^\omega$, the form $g_A$ is well-defined and nondegenerate, and
the Hamiltonian identities make $g_A$ and $b_L$ symmetric with
$g_A(C_Lx,y)=-b_L(x,y)$. Multiplying $(M-\lambda I)v=0$ on the left by
$I-\ii A$ yields $Bv+\lambda Av=\ii d_\lambda v$. Applying $A$ gives $ABv=\ii d_\lambda Av$. Applying $AB$ to the same equation and using $B^2=0$ yields $\lambda ABA v=\ii d_\lambda ABv=-d_\lambda^2 Av$, so $C_L(Av)=\rho Av$; conversely, any $x\in\ker(C_L-\rho I)$ lifts
to $v=(Bx+\ii d_\lambda x)/\rho$ with $Av=x$. Because $\rho\neq 0$, this
eigenspace lies in $\im C_L$, establishing \eqref{eq:pencil} and the rank bound.

Choosing a full-rank factorization $A=E_AF_A$ where $E_A$ is a basis matrix
for $L$, we have $F_AE_A=0$ and $F_ABE_A$ represents $C_L$. For
$d_\lambda\neq 0$, writing
$M-\lambda I=(I+\ii A)(d_\lambda I+\ii B+\ii\lambda A)$ and applying
$(d_\lambda I+\ii B)^{-1}=d_\lambda^{-1}I-\ii d_\lambda^{-2}B$ together with
$\det(I+XY)=\det(I+YX)$ yields the determinant identity, which extends to
$\lambda=1$ by continuity. Finally, at $\lambda=1$, the exact sequence
\eqref{eq:seq} and Remark~\ref{rem:forms} identify $\ker\mathcal O$ with
$\radform\bigl((g_A+g_B)|_{W_\cap\times W_\cap}\bigr)$, while
$\dim(\ker A\cap\ker B)=\dim(L+L_r)^\omega=2d_V-\rank A-\rank B+\dim W_\cap$.
\end{proof}

To compare Proposition~\ref{prop:shear} with Lemma~\ref{lem:ranktwo}, write
$B=E_BF_B^T$: then $C_L=(AE_B)(F_B^T|_L)$ whereas $C_1=(F_B^T|_L)(AE_B)$, so at
every $\rho\neq 0$ the maps $F_B^T|_L$ and $\rho^{-1}AE_B$ intertwine their
$\rho$-eigenspaces. Thus the real symmetric pencil tests the exact same
matching condition as the complex compression \eqref{eq:reduced}.

\begin{remark}[Lagrangian triple-form interpretation]\label{rem:triple}
For real $\rho$, consider the graphs of $b_L$ and $-\rho g_A$ in $L\oplus L^*$
together with the vertical Lagrangian $0\oplus L^*$ under the canonical pairing
$\omega_{\rm can}((x,p),(y,q))=q(x)-p(y)$. Subtracting $(x,-\rho g_A(x,\cdot))$
from $(x,b_L(x,\cdot))$ leaves the vertical covector $(b_L+\rho g_A)(x,\cdot)$,
so the Lagrangian triple form in this order is precisely $b_L+\rho g_A$ (in
particular, $b_L+4g_A$ at $\lambda=-1$). Whereas Remark~\ref{rem:fourier-graphs}
describes the Hermitian Schur complement for a physical twist, this real graph
form operates in the common real fiber of the two local monodromies; both
reductions detect the same Floquet kernel on their common domain, while
\eqref{eq:pencil} continues bilinearly to nonreal $\rho$.
\end{remark}

\subsection{Real-interval transport and proof of Theorem~B}\label{subsec:real-proof}

\begin{lemma}[Transport of the highest-exponent direction along $(r,0)$]\label{lem:real-transport}
For every real $\psi\in L$, let $x_\psi$ denote the position component of the
solution of \eqref{eq:fuchs-real} with state $\psi$ at $z_b$. Then $x_\psi(z)$
has a finite limit as $z\downarrow r$ along $(r,0)$, and
\begin{equation}\label{eq:sign-real}
 b_L(\psi,\psi)=\frac{2\pi}{\eta}
 \bigl(\alpha_2 x_{\psi,2}(r)^2-\alpha_1 x_{\psi,1}(r)^2\bigr).
\end{equation}
If $0<\alpha_1<\alpha_2$ and $\rank A=2$, then the vector $\phi\in L$ from
Lemma~\ref{lem:real-local} satisfies $g_A(\phi,\phi)=0$ and
$b_L(\phi,\phi)>0$.
\end{lemma}
\begin{proof}
In the local basis at $r$, a solution with coefficient vector $(p,x)^T$ has
position $x+O\bigl((z-r)|\log(z-r)|\bigr)$ because the lower block row of
$R_r$ in \eqref{eq:res-real} vanishes. Evaluating
$\omega\bigl((p,x),2\pi R_r(p,x)\bigr)$ and invoking invariance under real
symplectic transport along $(r,0)$ gives \eqref{eq:sign-real}.

Now assume $0<\alpha_1<\alpha_2$ and $\rank A=2$, and write
$x=x_\phi=(x_1,x_2)^T$. Set $\xi=-z$, $\xi_*=-r$, and
$f_e(\xi)=-1/\q(-\xi)>0$ for $0<\xi<\xi_*$. Denoting $\log\xi$-derivatives by
dots, the position equation becomes $\ddot x=2X\dot x+f_e(\xi)Dx$. By
Lemma~\ref{lem:real-local}, all four components of $(x,\dot x)$ are strictly
positive near $\xi=0$. Because the associated first-order system in $\log\xi$
has nonnegative off-diagonal entries, variation of constants keeps
$x_1,x_2,\dot x_1,\dot x_2>0$ on every compact subinterval of $(0,\xi_*)$.

Put $v_x=x_2-x_1$ and $\Delta_\alpha=\alpha_2-\alpha_1>0$. Subtracting the two
equations and switching to $\xi$-derivatives (denoted by primes) yields
$(\xi^3 v_x')'=\xi f_e(\xi)(\alpha_2 v_x+\Delta_\alpha x_1)$. Since
$v_x(\xi)=O(\xi^3)$ as $\xi\downarrow 0$, two integrations give the Volterra
equation
\[
 v_x(\xi)=\frac12\int_0^\xi\left(1-\frac{y^2}{\xi^2}\right)
 \frac{f_e(y)}{y}\bigl(\alpha_2 v_x(y)+\Delta_\alpha x_1(y)\bigr)\dd y.
\]
On each $[0,d]$ with $d<\xi_*$, the ratio $f_e(y)/y$ is bounded, so the $k$th
iterate of the Volterra operator has norm at most $(Cd)^k/k!$. Its convergent
Neumann series therefore sums nonnegative terms with a strictly positive
leading inhomogeneous term for $\xi>0$. Monotonicity of $x_1$ and existence
of the endpoint limits give $x_2(r)\ge x_1(r)>0$, whence $b_L(\phi,\phi)>0$ by
\eqref{eq:sign-real}, while $g_A(\phi,\phi)=0$ holds directly from
Lemma~\ref{lem:real-local}.
\end{proof}

\begin{lemma}[Positivity for scalar and small antiperiodic potentials]\label{lem:scalar}
Both antiperiodic vanishing assertions of \hyperref[main:real]{Theorem~B}\textup{(2)} hold. In
particular, $L_{I,e}$ has a strictly positive antiperiodic quadratic form for
every $0\le e<1$.
\end{lemma}
\begin{proof}
The $2\pi$-periodic rotation $v=\exp(tJ)u$ conjugates $L_{D,e}$ to
\[
 -\partial_t^2-I+\q^{-1}\exp(tJ)D\exp(-tJ)
\]
while preserving the antiperiodic domain. Integrating by parts and using
$\q''=1-\q$ gives the Picone-type identity
\begin{equation}\label{eq:picone}
 \int_0^{2\pi}\bigl(|v'|^2-|v|^2+|v|^2/\q\bigr)\dd t
 =\int_0^{2\pi}\q^2|(v/\q)'|^2\dd t.
\end{equation}
The right-hand side vanishes only when $v/\q$ is constant, which forces $v=0$
under antiperiodic boundary conditions; thus the form is strictly positive for
every scalar stiffness $D=\alpha I$ with $\alpha\ge 1$. On the antiperiodic
domain, half-integer Fourier frequencies also yield
$\norm{v'}_{L^2}^2\ge\frac14\norm{v}_{L^2}^2$, and since $\q<1+e$ almost
everywhere, strict positivity follows whenever $\alpha\ge 3(1+e)/4$. Finally,
the spectrum of $-\partial_t^2-I$ on the antiperiodic domain has distance $3/4$
from $0$. If $\norm{D}_{\rm op}\le 3(1-e)/4$, then $\q>1-e$ almost everywhere
gives $\norm{\q^{-1}Du}_{L^2}<\frac34\norm{u}_{L^2}$ for every nonzero $u$,
ruling out any nontrivial antiperiodic kernel.
\end{proof}

\begin{proof}[Proof of {\hyperref[main:real]{Theorem~B}}]
If $\min(\rank A,\rank B)\le 1$, then $\nu_\lambda(M)\le 1$ for every
$\lambda\neq 1$ follows immediately from Proposition~\ref{prop:shear}. Assume
therefore that $\rank A=\rank B=2$, so $\dim L=2$ and $g_A$ has signature
$(1,1)$ by Lemma~\ref{lem:real-local}. A two-dimensional eigenspace at
$\lambda\neq 1$ would force $b_L=-\rho g_A$ on all of $L$, with $\rho$ real and
nonzero because both forms are real and $g_A$ is nondegenerate. When $D$ is
indefinite, \eqref{eq:sign-real} makes $b_L$ positive semidefinite, which is
incompatible with signature $(1,1)$; when $D$ is positive definite and
nonscalar, Lemma~\ref{lem:real-transport} supplies a $g_A$-isotropic vector
$\phi$ with $b_L(\phi,\phi)>0$, again excluding proportionality.

For indefinite $D$, let $G$ and $H_b\ge 0$ represent $g_A$ and $b_L$ in a real
basis of $L$. By the Sylvester determinant identity, the nonzero eigenvalues of
$C_L=-G^{-1}H_b$ coincide with those of the real symmetric matrix
$-H_b^{1/2}G^{-1}H_b^{1/2}$ (and the argument holds vacuously when $\dim L=0$).
All values of $\rho$ in \eqref{eq:pencil} are therefore real, so
$\lambda+\lambda^{-1}=2-\rho$ places every multiplier in $\RR\cup\SSone$. Part~(2)
follows from Lemma~\ref{lem:scalar}, and since any real symmetric $2\times 2$
matrix with $\tr D\ge 2$ is either indefinite, singular, positive definite and
nonscalar, or scalar with $\alpha\ge 1$, the final assertion follows.
\end{proof}

\subsection{Three-body stiffnesses and the Euler--Lagrange antiperiodic curves}\label{subsec:three-body}

We now derive the shape stiffnesses of the planar three-body relative
equilibria in the exact normalization used later for regular polygons. Let
$m_i>0$ with Newton's gravitational constant set to one, and write the
gravitational acceleration vector $\mathcal F$ with components
\[
 \mathcal F_i(x)=\sum_{j\neq i}m_j\frac{x_j-x_i}{|x_j-x_i|^3}.
\]
Let $\mathbf a$ be a mass-centered central configuration satisfying
$\mathcal F(\mathbf a)=-\mu_\ast\mathbf a$ with $\mu_\ast>0$. Identifying the
physical plane with $\CC$, write $x(\tau)=R(\tau)e^{\ii t(\tau)}y(t(\tau))$ in
terms of physical time $\tau$ and true anomaly $t$, where $R=p_0/\q$ and
$R^2\dot t=\sqrt{\mu_\ast p_0}$ for a semilatus rectum $p_0>0$. Differentiating
gives $2\dot R\dot t+R\ddot t=0$ and $\ddot R-R\dot t^2=-\mu_\ast/R^2$, so
homogeneity and rotational covariance of $\mathcal F$ yield
\[
 R^3\dot t^2(y''+2Jy')=\mathcal F(y)+\mu_\ast y,\qquad R^3\dot t^2=\mu_\ast\q.
\]
Since the periodic state transformation preserves Floquet multiplicities, the
linearized position equation at $y=\mathbf a$ is
\begin{equation}\label{eq:newton-normalization}
 -u''-2Ju'+\q^{-1}(I+\mu_\ast^{-1}K)u=0,\qquad K=D\mathcal F(\mathbf a),
\end{equation}
where $J$ acts on each planar position component. The differential $K$ is
self-adjoint with respect to the mass inner product, satisfies
$K\mathbf a=2\mu_\ast\mathbf a$ by homogeneity and
$KJ\mathbf a=-\mu_\ast J\mathbf a$ by rotational invariance, and therefore
leaves invariant the orthogonal complement of the Kepler plane in the
center-of-mass space. For three bodies this essential shape complement is a
single real plane.

\begin{lemma}[Stiffness matrices of the collinear and equilateral configurations]\label{lem:three-shapes}
In a positively oriented orthonormal frame on the shape plane, the stiffness
matrix in \eqref{eq:newton-normalization} is $\diag(-\delta,3+2\delta)$ with
$\delta>0$ for a collinear three-body configuration, and $\diag(\kappa,3-\kappa)$
for the equilateral configuration, where
\[
 \beta=\frac{27(m_1m_2+m_2m_3+m_3m_1)}{(m_1+m_2+m_3)^2},\qquad
 \kappa=\frac{3-\sqrt{9-\beta}}{2}.
\]
Positive masses give $0<\beta\le 9$, and the equation extends analytically to
$\beta=0$.
\end{lemma}
\begin{proof}
For a collinear configuration ordered as $a_1<a_2<a_3$, put
$w_{ij}=|a_i-a_j|^{-3}$. On longitudinal and transverse displacements, $K$ acts
by $2K_{\rm col}$ and $-K_{\rm col}$, respectively, where
$(K_{\rm col}v)_i=\sum_{j\neq i}m_j w_{ij}(v_i-v_j)$. This operator annihilates
the constant vector and satisfies $K_{\rm col}\mathbf a=\mu_\ast\mathbf a$. In
the adjacent-gap coordinates $(v_2-v_1,v_3-v_2)$, the induced $2\times 2$
matrix $G$ has strictly negative off-diagonal entries
$G_{12}=m_3(w_{13}-w_{23})$ and $G_{21}=m_1(w_{13}-w_{12})$, while the
positive gap vector $g=(a_2-a_1,a_3-a_2)$ satisfies $Gg=\mu_\ast g$. Taking the
trace gives the second eigenvalue $k_\ast$ of $G$ via
\[
 k_\ast-\mu_\ast=-G_{12}\frac{g_2}{g_1}-G_{21}\frac{g_1}{g_2}>0.
\]
Ordering the transverse axis first on the shape plane,
\eqref{eq:newton-normalization} gives $\diag(-\delta,3+2\delta)$ with
$\delta=k_\ast/\mu_\ast-1>0$.

For the equilateral triangle normalized to side length one, put
$m_\Sigma=\sum_{i=1}^3 m_i$, so $\mu_\ast=m_\Sigma$. For a pair with unit
separation $u_{ij}=a_i-a_j$, differentiating $-|x|^{-3}x$ sends a relative
displacement $v$ to $(v+3u_{ij}^2\bar v)/2$. In mass-weighted coordinates
$z_i=\sqrt{m_i}\,u_i$, the stiffness operator $I+\mu_\ast^{-1}K$ on the
center-of-mass space therefore acts by
\[
 z\longmapsto\tfrac32 z+\tfrac{3}{2m_\Sigma}B_\triangle\bar z,\qquad
 (B_\triangle)_{ii}=\sum_{j\neq i}m_j u_{ij}^2,\quad
 (B_\triangle)_{ij}=-\sqrt{m_i m_j}\,u_{ij}^2\quad(i\neq j).
\]
The complex symmetric matrix $B_\triangle$ defines a real self-adjoint
antilinear map that annihilates translations and sends the mass-weighted
configuration vector to $m_\Sigma$ times itself, so its singular value on the
Kepler complex line is $m_\Sigma$. On the orthogonal complex shape line it acts
in a unit coordinate by $z\mapsto c_\ast\bar z$, with real eigenvalues
$\pm|c_\ast|$. Since the squared unit directions at a common vertex have real
inner product $\cos(2\pi/3)=-1/2$, computing the Frobenius norm squared yields
\[
 m_\Sigma^2+|c_\ast|^2=\norm{B_\triangle}_{\rm F}^2
 =2\sum_{i=1}^3 m_i^2+\sum_{i<j}m_i m_j,
\]
whence $|c_\ast|^2=m_\Sigma^2-3\sum_{i<j}m_i m_j$ and the two stiffness
eigenvalues are $\frac32\bigl(1\pm\sqrt{1-\beta/9}\bigr)$, including the
equal-mass case $c_\ast=0$.
\end{proof}

\begin{lemma}[Invertibility and Morse index at the Kepler endpoint]\label{lem:ap-kepler-endpoint}
For $D=\diag(0,3)$ and every $0\le e<1$, the antiperiodic operator $L_{D,e}$ is
invertible and has Morse index two.
\end{lemma}
\begin{proof}
Writing $u=(x,y)^T$, the first equation of $L_{D,e}u=0$ makes $x'-2y$ constant,
hence zero under antiperiodic boundary conditions, and the second reduces to
$y''+(4-3/\q)y=0$. Because $(\partial_t^2+4)(\q\sin t)=3\sin t$, the nonzero
function $\q(t)\sin t$ is a $2\pi$-periodic solution. Since this scalar
equation has constant Wronskian, its monodromy has determinant one and
eigenvalue $1$, hence both eigenvalues equal $1$; consequently $y=0$ for every
antiperiodic solution, which forces $x'=0$ and thus $x=0$.

At $e=0$, the complex Fourier blocks at half-integer frequencies
$k\in\mathbb Z+\frac12$ are $k^2 I-2\ii kJ+D$, with positive trace and
determinant $k^2(k^2-1)$. Only $k=\pm\frac12$ contribute a negative eigenvalue
(one from each block), so the Morse index is two at $e=0$. On the fixed
antiperiodic domain, the bounded potential varies norm-continuously in
$e\in[0,1)$ without ever creating a zero eigenvalue, so the Morse index remains
two for every $0\le e<1$.
\end{proof}

\begin{proof}[Proof of {\hyperref[main:real-curves]{Corollary~B.1}}]
Consider an affine family $D=D_0+sS$ where $D_0$ is a positive scalar multiple
of $I$ such that $R_e=L_{D_0,e}$ is strictly positive on the antiperiodic
domain, and put $V_e=\q^{-1}S$ and $K_e=R_e^{-1/2}V_eR_e^{-1/2}$. Because
$R_e$ has compact resolvent, $K_e$ is compact and self-adjoint, and the
isomorphism $u\mapsto R_e^{1/2}u$ from the antiperiodic form domain onto $L^2$
sends the quadratic form of $L_{D,e}$ to that of $I+sK_e$. The positive values
of $s$ for which $\ker L_{D,e}\neq\{0\}$, counted with multiplicity, are
therefore $s_j(e)=-1/\mu_j(e)$, where $\mu_1(e)\le\mu_2(e)\le\dots<0$ are the
negative eigenvalues of $K_e$, and the Morse index at $s>0$ counts the values
$s_j(e)<s$. In both families below, $S$ is diagonal with negative first entry,
so smooth antiperiodic test functions supported in the first coordinate span
arbitrarily large negative subspaces for $V_e$ and hence for $K_e$, ensuring an
infinite sequence $s_j(e)$. By the resolvent identity and continuity of the
positive square root, $K_e$ depends norm-continuously on $e\in[0,1)$, so each
ordered root $s_j(e)$ is continuous on $[0,1)$ by the min--max principle.

For the extended Euler family, choose $D_0=I$, $S=\diag(-1,2)$, and
$s=1+\delta$, with $R_e>0$ by Lemma~\ref{lem:scalar}. At $s=1$ we obtain
$D=\diag(0,3)$, so Lemma~\ref{lem:ap-kepler-endpoint} places exactly two roots
below $1$ and none at $1$; thus the Euler roots in $\delta>0$ are precisely
$s_{j+2}(e)-1$ for $j\ge 1$. At $0<e<1$, \hyperref[main:real]{Theorem~B} forces every root to be
simple, making the sequence strictly increasing. At $e=0$, the $k$th
half-integer Fourier block has determinant $k^4+(\delta-1)k^2-2\delta^2-3\delta$.
For $k=\pm\frac12$ this is strictly negative for all $\delta\ge 0$, whereas for
each pair $k=\pm(j+\frac12)$ ($j\ge 1$) it has a unique positive root that
increases strictly with $j$ because the determinant increases with $k^2>1$ at
fixed $\delta\ge 0$. Each circular value therefore has multiplicity two, so
$s_{2j+1}(e)-1<s_{2j+2}(e)-1$ are the two strictly separated branches emerging
from the $j$th circular pair.

For the Lagrange family, take $D_0=\frac32 I$, $S=\diag(-1,1)$, and
$s=\sqrt{9-\beta}/2\in[0,3/2]$, where again $R_e>0$ by Lemma~\ref{lem:scalar}
(so in particular $\beta=9$, corresponding to $s=0$, carries no kernel). At
$s=3/2$ ($\beta=0$), we recover the Kepler endpoint $\diag(0,3)$, so
Lemma~\ref{lem:ap-kepler-endpoint} yields exactly two roots in $(0,3/2)$ and
hence exactly two Lagrange parameter values in $(0,9)$, ordered inversely via
$\beta=9-4s^2$. For $0<e<1$ they are distinct by \hyperref[main:real]{Theorem~B} because $\tr D=3$.
At $e=0$, the Fourier determinant $k^4-k^2+\beta/4$ vanishes only at
$k=\pm\frac12$ and $\beta=3/4$, where both branches coalesce.

Finally, any nonzero kernel vector $u$ at $s>0$ satisfies
$\langle u,V_eu\rangle=-s^{-1}\langle u,R_eu\rangle<0$. At $0<e<1$ the kernel
is one-dimensional, so taking the Schur complement on $u^\perp$ reduces the
analytic fixed-domain family to a scalar analytic equation with nonvanishing
$s$-derivative. The implicit function theorem yields local real-analytic
branches $s(e)$ (with negative crossing sign in $\delta$ and positive crossing
sign in $\beta$), which glue globally along the continuous ordered enumeration.
\end{proof}

\section{Complex gyroscopic blocks and the trace difference}\label{sec:complex}

We turn now to the complex two-component equation \eqref{eq:complex-op}.
Conjugating by $\diag(1,-1)$ reverses the sign of $c$ while leaving $Z$
invariant, so throughout the proofs in this section we may assume $c>0$. From
the viewpoint of the physical eight-dimensional real symplectic system of a
polygonal mode, the complexification in $\mathrm{Sp}(8,\CC)$ splits into
\eqref{eq:complex-op} paired with its transpose-dual equation. Two distinct
transports enter our analysis: along the real segment $(r,0)$ both equations
have real coefficients and distinguished real positive solutions, whereas along
real physical time $t\in\RR$ complex conjugation identifies the two solution
spaces and produces a conserved indefinite Hermitian form.

\begin{proposition}[Conserved Hermitian pairing on the physical circle]\label{prop:current}
Consider the direct equation $w''-2\ii Zw'=\q^{-1}Hw$ and its transpose-dual
equation $v''+2\ii Zv'=\q^{-1}Hv$. Their bilinear concomitant
\[
 W_t(v,w)=v^Tw'-(v')^Tw-2\ii v^TZw
\]
is independent of $t$. Writing $Y_w=(w,-\ii w')^T$ along real time, any two
solutions $u,w$ of \eqref{eq:complex-op} preserve the Hermitian form
(specializing \eqref{eq:common-current} to $s=2$)
\begin{equation}\label{eq:current-matrix}
 [u,w]_{\mathcal J}=-\ii W_t(\overline u,w)=Y_u^*\mathcal J Y_w,\qquad
 \mathcal J=\begin{pmatrix}-2Z&I\\I&0\end{pmatrix}.
\end{equation}
In physical state coordinates, $\mathcal J$ is nondegenerate of signature
$(2,2)$, and the period map satisfies $M^*\mathcal J M=\mathcal J$,
$\det M=1$, and $(\im(M-I))^{\perp_{\mathcal J}}=\ker(M-I)$.
\end{proposition}
\begin{proof}
Differentiating $W_t(v,w)$ cancels the $H$-terms by symmetry of $H$, while the
drift terms combine into $2\ii(v^TZw)'$. Along real $t$, $\overline u$ solves
the transpose-dual equation whenever $u$ solves \eqref{eq:complex-op}, yielding
\eqref{eq:current-matrix}. The matrix $\mathcal J$ is nondegenerate and has
signature $(2,2)$ because the last two standard basis vectors span an isotropic
two-plane. Since the first-order time system has trace zero, Liouville's
formula gives $\det M=1$. Finally, $M^*\mathcal J M=\mathcal J$ identifies the
$\mathcal J$-adjoint of $M$ with $M^{-1}$, whence
$(\im(M-I))^{\perp_{\mathcal J}}=\ker(M^{-1}-I)=\ker(M-I)$.
\end{proof}

\begin{lemma}[Origin resonances and nilpotent generators]\label{lem:resonance}
In the real $z$-frame $Y=(w,\vang w)^T$, the based local factors of
\eqref{eq:complex-op} at $z_b\in(r,0)$ take the form $U_0=\exp(\ii A)$ and
$U_r=I+\ii B$, where $A,B$ are real matrices satisfying $A^3=0$, $B^2=0$, and
$\rank B=\rank H$. Before transport to $z_b$ and multiplication by $2\pi$, the
origin logarithm in the exponent basis ordered as $(-2,0,0,2)$ has the form
\begin{equation}\label{eq:generator}
 N=\begin{pmatrix}
  0&0&0&0\\ \ru_1&0&0&0\\ \ru_2&0&0&0\\ *&\rv_1&\rv_2&0
 \end{pmatrix},\qquad
 \begin{aligned}
  \ru_1&=-\frac{2c(a-3b+3)}{3e^2},&\quad
  \ru_2&=-\frac{2(3b^2-3b-c^2)}{3e^2},\\
  \rv_1&=-\frac{2(3a^2-3a-c^2)}{3e^2},&\quad
  \rv_2&=-\frac{2c(3a-b-3)}{3e^2}.
 \end{aligned}
\end{equation}
In particular, $N^3=0$ and
\begin{equation}\label{eq:K}
 N^2=\mathcal K E_{41},\qquad
 \mathcal K=\rv_1\ru_1+\rv_2\ru_2=\frac{4c(a-b)}{9e^4}\mathcal T(H),
\end{equation}
where $E_{41}$ is the $(4,1)$ matrix unit and $\mathcal T(H)$ is given by
\eqref{eq:T}.
\end{lemma}
\begin{proof}
In the complex time variable $z=e^{\ii t}$, equation \eqref{eq:complex-op} and
its residue at the inner pole $r$ are
\[
 \vang^2w-2Z\vang w+\frac{2z}{ez^2+2z+e}Hw=0,\qquad
 R_r=\begin{pmatrix}0&0\\-H/\eta&0\end{pmatrix}.
\]
Because $R_r^2=0$, every positive homological order at $r$ is invertible by
Lemma~\ref{lem:gauge-rec}, so the pole logarithm is conjugate to $R_r$ and
hence $B^2=0$ with $\rank B=\rank H$. At the origin, choose exponent
eigenvectors with leading positions
$(\mathbf e_2,\mathbf e_1,\mathbf e_2,\mathbf e_1)$ in exponent order
$(-2,0,0,2)$ and lower component equal to the exponent times the position.
Positive exponent differences are $2$ and $4$, with only the $(4,1)$ entry
resonant at order four. Appendix~\ref{app:gauge} proves convergence of the
gauge and evaluates the order-two resonance coefficients
$\ru_1,\ru_2,\rv_1,\rv_2$ in \eqref{eq:generator}; direct multiplication of the
triangular matrix $N$ gives \eqref{eq:K}. Since the equation has real Taylor
coefficients on $(r,0)$, choosing real logarithm branches on the two stems
yields real local fundamental matrices, whose values at $z_b$ conjugate $2\pi N$
and $2\pi R_r$ to real matrices $A$ and $B$.
\end{proof}

Because $A^3=0$ and $B^2=0$, the exact Pad\'e identities
$\exp(\ii A)=(I+\ii A/2)(I-\ii A/2)^{-1}$ and
$I+\ii B=(I+\ii B/2)(I-\ii B/2)^{-1}$ identify $A$ and $B$ with the Cayley
generators of Theorem~\ref{thm:cayley}. It remains to express the trace difference $\Delta(M)$
in terms of transport along $(r,0)$.

\begin{lemma}[Positive transport of the dual highest-exponent solutions]\label{lem:dual-transport}
Assume $a,b,c>0$, and let $w_+(z)=z^2(\mathbf e_1+O(z))$ and
$v_+(z)=z^2(\mathbf e_2+O(z))$ be the unique highest-exponent solutions of the
direct and transpose-dual equations at the origin. Both have strictly positive
position components at $r$, and for each fixed $0<e<1$ the connection scalar
\[
 F(H,e)=\frac{v_+(r)^THw_+(r)}{\eta}
\]
is strictly positive and depends continuously on the entries of $H$ in the
region $a,b,c>0$.
\end{lemma}
\begin{proof}
Since exponent $2$ has no positive resonance, $w_+$ and $v_+$ are uniquely
determined. Put $\xi=-z$ and $k_e(\xi)=2/[\xi(e\xi^2-2\xi+e)]>0$ on
$(0,\xi_*)$, where $\xi_*=-r$. Writing primes for $\xi$-derivatives, the
components of $w_+$ obey
\[
 (w_1'/\xi)'=\frac{k_e(\xi)}{\xi}(aw_1+cw_2),\qquad
 (\xi^3w_2')'=\xi^3 k_e(\xi)(cw_1+bw_2).
\]
Direct series substitution yields $w_1(\xi)=\xi^2+O(\xi^3)$ and
$w_2(\xi)=\frac{2c}{15e}\xi^3+O(\xi^4)$, so $w_1,w_2,w_1',w_2'$ are initially
positive near $\xi=0$. Integrating the two weighted equations prevents any
first zero of either component or derivative on $(0,\xi_*)$, and interchanging
the two weights and components gives the same positivity for $v_+$. Because the
position block row of $R_r$ vanishes, both increasing positions possess finite
strictly positive limits at $r$; since all entries of $H$ are positive,
$F(H,e)>0$ holds even when $H$ is indefinite.

On any compact neighborhood in $\{a,b,c>0\}$, both the highest-exponent Taylor
recursion at $0$ and the nonresonant gauge recursion at $r$ admit uniform
geometric majorants by Lemma~\ref{lem:gauge-rec}. Propagation between two fixed
ordinary points in $(r,0)$ depends continuously on $H$, and in the local pole
basis the endpoint limit $w_+(r)$ is simply the position part of the
coefficient vector, yielding continuity of $F(H,e)$ in $H$.
\end{proof}

\begin{proposition}[Trace difference and the connection pairing]\label{prop:trace}
Under the hypotheses of Lemma~\ref{lem:dual-transport}, the period map $M$
satisfies
\begin{equation}\label{eq:imtrace}
 \tr(A^2B)=4\pi^3\mathcal K F(H,e),\qquad
 \Delta(M)=2\ii\imag\tr M=-4\pi^3\ii\mathcal K F(H,e).
\end{equation}
\end{proposition}
\begin{proof}
Analytically continuing $W_t/\ii$ from Proposition~\ref{prop:current} produces
the constant bilinear pairing
\[
 W_z(v,w)=v^T\vang w-(\vang v)^Tw-2v^TZw.
\]
Let $\ell_-(w)$ extract the coefficient of the exponent $-2$ basis solution in
the origin Frobenius frame of the direct equation. Evaluating $W_z(v_+,w)$ on
leading terms at $z=0$ pairs $v_+$ with $z^{-2}\mathbf e_2$ to $-2$ and
annihilates the remaining columns (including logarithmic terms) by their orders
at the origin, so $\ell_-(w)=-W_z(v_+,w)/2$. In the common fiber at $z_b$,
\eqref{eq:K} therefore gives $A^2=(2\pi)^2\mathcal K\,w_+\otimes\ell_-$, where
$w_+$ denotes the state of $w_+$ at $z_b$.

Because the pole gauge equals $I$ at $r$ and $R_r^2=0$, the solution $Bw_+$ has
position of order $O(z-r)$ at $r$: the Frobenius vector of $Bw_+$ lies in $\im R_r\subset\ker R_r$, so its logarithmic term vanishes. Its lower
component tends to $-2\pi Hw_+(r)/\eta$. In $W_z(v_+,Bw_+)$, the product of
the logarithmically growing derivative of $v_+$ with the vanishing position of
$Bw_+$ tends to zero as $z\downarrow r$, leaving
$\ell_-(Bw_+)=-W_z(v_+,Bw_+)/2=\pi F(H,e)$. Taking the trace of the rank-one
operator $A^2B$ yields the first identity in \eqref{eq:imtrace}. Theorem~\ref{thm:cayley}(2)
gives $\Delta(M)=-\ii\tr(A^2B)$, while $M^*\mathcal J M=\mathcal J$ in
Proposition~\ref{prop:current} implies $\tr M^{-1}=\overline{\tr M}$, whence
$\Delta(M)=2\ii\imag\tr M$. Since traces are invariant under similarity, these
identities hold equally in the physical and based frames.
\end{proof}

\begin{lemma}[Transverse stiffness variation on $\mathcal T(H)=0$]\label{lem:tangent}
Under the hypotheses of \hyperref[main:complex]{Theorem~C}, if $\mathcal T(H)=0$ then
$\nu_\lambda(M)\le 1$ for every $\lambda\neq 1$.
\end{lemma}
\begin{proof}
Fix $0<e<1$ and consider the real family $H_\varepsilon=H+\varepsilon I$ for
small real $\varepsilon$, maintaining $c>0$ by the initial sign normalization.
By \eqref{eq:T} and \eqref{eq:K}, $\mathcal K(0)=0$ and
$\mathcal K'(0)=48c(a-b)/(9e^4)\neq 0$. Because $\mathcal K$ vanishes at
$\varepsilon=0$, continuity of $F(H_\varepsilon,e)$ from
Lemma~\ref{lem:dual-transport} suffices to differentiate \eqref{eq:imtrace} at
$\varepsilon=0$, giving
\[
 \left.\frac{\mathrm d}{\mathrm d\varepsilon}\Delta(M(\varepsilon))\right|_0
 =-4\pi^3\ii\mathcal K'(0)F(H,e)\neq 0.
\]
Along this smooth physical family, Proposition~\ref{prop:current} guarantees
$\det M(\varepsilon)=1$, and at $\varepsilon=0$ equation \eqref{eq:K} gives
$A^2=0$. Applying Theorem~\ref{thm:cayley}(2b) to the smooth physical period map $M(\varepsilon)$
rules out any two-dimensional eigenspace at $\lambda\neq 1$.
\end{proof}

\begin{proof}[Proof of {\hyperref[main:complex]{Theorem~C}}]
After normalizing $c>0$, if $\mathcal T(H)\neq 0$ then \eqref{eq:K} and
Proposition~\ref{prop:trace} give $\Delta(M)\neq 0$, so Theorem~\ref{thm:cayley}(2a) yields
both $\nu_\lambda(M)\le 1$ for all $\lambda\neq 1$ and algebraic simplicity at
$\lambda=\pm 1$. When $\mathcal T(H)=0$, Lemma~\ref{lem:tangent} establishes
$\nu_\lambda(M)\le 1$ for every $\lambda\neq 1$.
\end{proof}

\begin{remark}[Obstruction to single-block complex symplecticity]\label{rem:no-single-symp}
In complex time, physical conjugation corresponds to inversion $z\mapsto 1/\bar z$
fixing the unit circle pointwise, whereas reflection across the real connecting
segment $(r,0)$ is $z\mapsto\bar z$. Whenever $\imag\tr M\neq 0$, the
characteristic polynomial of the single four-dimensional block $M$ is not
reciprocal, so $M$ cannot preserve any nondegenerate complex alternating
bilinear form on this block alone; rather, the ambient complexified symplectic
structure lives on the eight-dimensional direct sum of the block and its
transpose dual, where it restricts to the bilinear concomitant $W_z(v,w)$.
\end{remark}

\section{Physical Fourier blocks of regular \texorpdfstring{$(1+n)$}{(1+n)}-gons}\label{sec:polygon}

Consider a planar regular $n$-gon of unit ring masses together with a central
mass $m>0$. Define the polygonal trigonometric sums
\begin{equation}\label{eq:sums}
\begin{gathered}
 \theta_j=\frac{2\pi j}{n},\qquad d_j=2\sin\frac{\pi j}{n},\qquad
 \gamma_n=\sum_{j=1}^{n-1}\frac{1}{2d_j},\qquad \mu=m+\gamma_n,\\
 P_\ell=\sum_{j=1}^{n-1}\frac{1-\cos(\ell\theta_j)\cos\theta_j}{2d_j^3},\quad
 Q_\ell=\sum_{j=1}^{n-1}\frac{\cos\theta_j-\cos(\ell\theta_j)}{2d_j^3},\quad
 S_\ell=\sum_{j=1}^{n-1}\frac{\sin(\ell\theta_j)\sin\theta_j}{2d_j^3}.
\end{gathered}
\end{equation}
We derive the physical Fourier blocks of Hu, Long, and Ou \cite{HLO} directly
from \eqref{eq:newton-normalization}. Identify the plane with $\CC$ so that $J$
acts as multiplication by $\ii$, set $\omega_n=e^{2\pi\ii/n}$, and take the
reference configuration $a_{\rm c}=0$, $a_j=\omega_n^j$ ($0\le j<n$). At $a_j$,
the ring acceleration divided by $\omega_n^j$ has zero imaginary part and real
part $\sum_{r=1}^{n-1}(\cos\theta_r-1)/d_r^3=-\gamma_n$, so
$\mathcal F(\mathbf a)=-\mu\mathbf a$, including at the center by symmetry.
Equation \eqref{eq:newton-normalization} therefore gives the stiffness operator
$I+\mu^{-1}K$, where $K=D\mathcal F(\mathbf a)$.

For ring displacements $u_j=\omega_n^j\xi_j$ and central displacement
$u_{\rm c}$, pass to the mass-unitary Fourier coordinates
\[
 \xi_j=n^{-1/2}\sum_{k\bmod n}X_k\omega_n^{kj},\qquad W=\sqrt{m}\,u_{\rm c},
\]
so that $|W|^2+\sum_k|X_k|^2$ is the physical mass norm squared; here each
$X_k$ is a complex planar coordinate, independent of $X_{-k}$ whenever
$k\not\equiv -k\pmod n$. A pair separated by $d e^{\ii\phi}$ contributes
$(v+3e^{2\ii\phi}\bar v)/(2d^3)$ to the inverse-distance Hessian on a relative
displacement $v$. Between ring vertices $j$ and $j+s$, the local phase factor
is $(1-\omega_n^s)^2/d_s^2=-\omega_n^s$, so the ring contribution to $(Ku)_j$
in the local frame at $a_j$ is
\[
 \sum_{s=1}^{n-1}\frac{\xi_j-\omega_n^s\xi_{j+s}-3\omega_n^s\bar\xi_j+3\bar\xi_{j+s}}{2d_s^3},
\]
while the central mass contributes
$\frac{m}{2}(\xi_j+3\bar\xi_j-\omega_n^{-j}u_{\rm c}-3\omega_n^j\bar u_{\rm c})$.
Taking Fourier coefficients (and using $\sum_{j=0}^{n-1}\omega_n^{2j}=0$ for
$n\ge 3$) yields
\[
\begin{aligned}
 (KX)_k&=A_kX_k+C_k\bar X_{-k}
 -\frac{\sqrt{mn}}{2}\delta_{k,-1}W-\frac{3\sqrt{mn}}{2}\delta_{k,1}\bar W,\\
 KW&=\frac{n}{2}W-\frac{\sqrt{mn}}{2}X_{-1}-\frac{3\sqrt{mn}}{2}\bar X_1,\\
 A_k&=\frac{m}{2}+P_k+S_k,\qquad C_k=\frac32(m-2Q_k).
\end{aligned}
\]
The mode $k=0$ is the Kepler plane, where $K$ acts diagonally with eigenvalues
$2\mu,-\mu$ and normalized stiffness $\diag(3,0)$.

In the first mode ($k=\pm 1$), split $(X_{-1},W)$ into the translation
coordinate $(\sqrt{n}\,X_{-1}+\sqrt{m}\,W)/\sqrt{m+n}$ (annihilated by $K$) and
the orthogonal shape coordinate
$Y=(\sqrt{m}\,X_{-1}-\sqrt{n}\,W)/\sqrt{m+n}$. Since $Q_1=0$ and $S_1=P_1$,
the remaining two equations read
\[
 KY=\frac{m+n}{2}Y+\frac32\sqrt{m(m+n)}\,\bar X_1,\qquad
 (KX)_1=\left(\frac{m}{2}+2P_1\right)X_1+\frac32\sqrt{m(m+n)}\,\bar Y.
\]
Setting $w=(\bar Y,X_1)^T$ gives \eqref{eq:complex-op} with gyroscopic term
$2\ii Z\partial_t$ and stiffness entries
\begin{equation}\label{eq:first-mode}
 a=1+\frac{m+n}{2\mu},\qquad b=1+\frac{m+p_n}{2\mu},\qquad
 c=\frac{3\sqrt{m(m+n)}}{2\mu},\qquad p_n=4P_1.
\end{equation}
For each higher mode $2\le\ell<n/2$, setting $w=(\bar X_{-\ell},X_\ell)^T$
again produces \eqref{eq:complex-op} with stiffness entries
\begin{equation}\label{eq:high-mode}
 a=\frac{3m/2+\gamma_n+P_\ell-S_\ell}{\mu},\qquad
 b=\frac{3m/2+\gamma_n+P_\ell+S_\ell}{\mu},\qquad
 c=\frac{3(m-2Q_\ell)}{2\mu}.
\end{equation}
When $n$ is even, the terminal index $\ell=n/2$ satisfies $\ell\equiv -\ell\pmod n$
and $S_\ell=0$, so $X_\ell$ spans a single real two-component position plane
governed by \eqref{eq:real-op} with diagonal stiffness
\[
 D_{n/2}=\diag\!\left(1+\frac{P_\ell-3Q_\ell+2m}{\mu},
 1+\frac{P_\ell+3Q_\ell-m}{\mu}\right).
\]
Together with their velocities, these coordinates decompose the full essential
state space (for $n=3$, the first complex mode accounts for the entire
eight-dimensional real essential state). Because each complex mode is the
realification of its four-dimensional complex state space, at every real
multiplier $\lambda$ one has
$\dim_{\RR}\ker(M_{\RR}-\lambda I)=2\dim_{\CC}\ker(M-\lambda I)$, whereas at
the real terminal block complexification preserves the numerical nullity.

\begin{lemma}[Signs and finite differences of the polygonal sums]\label{lem:polygon-sums}
For $3\le n\le 11$ one has $p_n<n$, and for $n\ge 12$ one has $p_n>n$;
moreover, $4\gamma_n+p_n>2n$ for every $n\ge 7$. For $n\ge 4$, the sum $Q_j$ is
strictly positive for every $2\le j\le\lfloor n/2\rfloor$. Finally, for every
$2\le\ell<n/2$ one has $S_\ell>0$ and
\begin{equation}\label{eq:polygon-addition}
 P_\ell\pm S_\ell=\frac{\gamma_n}{2}+Q_{\ell\pm 1},\qquad Q_{\ell\pm 1}\ge 0.
\end{equation}
For $t_n=\pi/(2n)$ and $2\le j\le\lfloor n/2\rfloor+1$ the moments satisfy
\begin{equation}\label{eq:cot-moments}
 Q_j=\frac{j^2-1}{2}\gamma_n
 -\frac14\sum_{i=0}^{j-2}(j-i)(j-i-1)\cot((2i+1)t_n).
\end{equation}
\end{lemma}
\begin{proof}
Since $1-\cos^2\theta_j=\sin^2\theta_j$, we may write
$p_n=\sum_{j=1}^{n-1}f(j/n)$ and $4\gamma_n=\sum_{j=1}^{n-1}g(j/n)$, where
$f(x)=\csc(\pi x)-\sin(\pi x)$ and $g(x)=\csc(\pi x)$. Both functions are
strictly convex on $(0,1)$, as
$g''(x)=\pi^2\csc(\pi x)\bigl(\cot^2(\pi x)+\csc^2(\pi x)\bigr)>0$ and
$f''(x)=g''(x)+\pi^2\sin(\pi x)>0$. Each grid point $j/n$ is the convex
combination $(1-j/n)\frac{j}{n+1}+(j/n)\frac{j+1}{n+1}$, so summing the strict
convexity inequalities for either grid sum $G_n$ yields
$G_n<\frac{n-1}{n}G_{n+1}$. Because $G_n>0$ for $n\ge 3$ and
$\frac{n-1}{n}<\frac{n}{n+1}$, the normalized sequence $G_n/n$ is strictly
increasing. It therefore suffices to verify $p_{11}<11$, $p_{12}>12$, and
$4\gamma_7+p_7>14$, which is done via elementary rational bounds in
Lemma~\ref{lem:finite-first-mode} of Appendix~\ref{app:finite-n}.

Put $M_n=\lfloor n/2\rfloor$, $h_j=Q_{j+1}-Q_j$ and
$B_j=\sum_{k=1}^{n-1}\cos(j\theta_k)/(2d_k)$. The cosine addition formula
and $d_k^2=2(1-\cos\theta_k)$ give \eqref{eq:polygon-addition} and
$h_j-h_{j-1}=B_j$. With $Q_0=-\gamma_n/2$, $Q_1=0$ and $B_0=\gamma_n$,
the first differences of $B_j$ are
\[
 B_j-B_{j+1}=\tfrac12\sum_{k=1}^{n-1}\sin((2j+1)\pi k/n)
            =\tfrac12\cot((2j+1)t_n)>0\qquad(0\le j<M_n).
\]
The sine sum telescopes after multiplication by $2\sin((2j+1)t_n)$.
Summing the two difference relations gives \eqref{eq:cot-moments}.

Since $h_0=B_0/2$, the quotient $h_j/(2j+1)$ is the weighted mean of
$B_0/2,\ldots,B_j/2$ with weights $1,2,\ldots,2$, and therefore decreases
strictly up to $j=M_n$. By reflection $Q_{n-j}=Q_j$, we have $h_{M_n}=0$ for
odd $n$ and $h_{M_n}=-h_{M_n-1}$ for even $n$, forcing $h_j>0$ for all
$0\le j<M_n$ in both cases. Consequently, $Q_j>0$ for $2\le j\le M_n$,
$2S_\ell=Q_{\ell+1}-Q_{\ell-1}=h_{\ell-1}+h_\ell>0$ for every $2\le\ell<n/2$
(including $\ell=M_n$ when $n$ is odd), and reflection also yields the endpoint
nonnegativity in \eqref{eq:polygon-addition}.
\end{proof}

\begin{proof}[Proof of {\hyperref[main:physical]{Corollary~C.1}}]
In the first complex mode \eqref{eq:first-mode}, Lemma~\ref{lem:polygon-sums}
gives $a,b>0$, $c>0$, and $a-b=(n-p_n)/(2\mu)\neq 0$ for every $n\ge 3$, so
\hyperref[main:complex]{Theorem~C} yields $\nu_\lambda(M)\le 1$ for all $m>0$, $0<e<1$, and
$\lambda\neq 1$. Expanding \eqref{eq:T} yields
\[
 \mu^2\mathcal T(H)=3\left[(4\gamma_n+p_n-2n)m+\gamma_n^2
 +(n+p_n)\gamma_n+\frac{(n-p_n)^2}{4}\right],
\]
which is strictly positive for every $n\ge 7$ by Lemma~\ref{lem:polygon-sums};
hence \hyperref[main:complex]{Theorem~C} implies algebraic simplicity at $+1$ and $-1$ throughout this
range.

In a higher complex mode $2\le\ell<n/2$, identity \eqref{eq:polygon-addition}
rewrites \eqref{eq:high-mode} as
\[
 a=\frac32+\frac{Q_{\ell-1}}{\mu},\qquad
 b=\frac32+\frac{Q_{\ell+1}}{\mu},\qquad
 a-b=-\frac{2S_\ell}{\mu}\neq 0,
\]
so \hyperref[main:complex]{Theorem~C} applies at every coupled mass $m\neq 2Q_\ell$. At the decoupling
mass $m=2Q_\ell$, the two components separate, and the $2\pi$-periodic gauges
$e^{\pm\ii t}$ remove the scalar drifts while preserving any unit twist
$\lambda\in\SSone$, reducing each component to $-\partial_t^2-1+A/\q$ with
$A=a$ or $b$. Because boundary terms cancel for every $|\lambda|=1$, the
identity \eqref{eq:picone} holds on every unit-twist domain; since
$A\ge 3/2>1$, it yields the strict form bound
\begin{equation}\label{eq:decoupled-pos}
 \int_0^{2\pi}\bigl(|v'|^2-|v|^2+A\q^{-1}|v|^2\bigr)\dd t
 \ge\frac{A-1}{1+e}\int_0^{2\pi}|v|^2\dd t>0.
\end{equation}
Thus the operator at $m=2Q_\ell$ is invertible for every $|\lambda|=1$, so in
particular every antiperiodic (and every periodic) kernel in a higher complex
mode occurs at a coupled mass, yielding $\nu_{-1}(M)\le 1$ for all $m>0$.

Finally, for even $n$ and $\ell=n/2$, the real terminal block has diagonal
stiffness with trace $2+(m+2P_\ell)/\mu>2$ because $P_\ell\ge 0$ termwise, so
\hyperref[main:real]{Theorem~B} yields $\nu_{-1}(M)\le 1$ for every $m>0$.
\end{proof}

\begin{remark}[Mass-crossing signs in higher modes]\label{rem:crossings}
The decoupling mass $m_s=2Q_\ell$ also fixes the direction of antiperiodic
crossings. From \eqref{eq:high-mode}, the rescaled family
$B_m=(m+\gamma_n)\mathcal L_{H(m),e}$ is affine in $m$ on a fixed domain. For
any nonzero $u$ in the antiperiodic kernel at $m$, evaluating
$B_{m_s}=B_m+(m_s-m)\partial_m B_m$ against $u$ yields
\[
 (m+\gamma_n)(m_s-m)\langle u,\partial_m\mathcal L_{H(m),e}u\rangle
 =(m_s+\gamma_n)\langle u,\mathcal L_{H(m_s),e}u\rangle>0
\]
by \eqref{eq:decoupled-pos}. Consequently, every simple antiperiodic root
generates a locally real-analytic mass branch $m(e)$ with opposite crossing
signs on the two sides of $m_s$.
\end{remark}

\section{Periodic eigenspaces on the resonant surface \texorpdfstring{$\mathcal T=0$}{T=0}}\label{sec:periodic}

At $\lambda=1$, the compression in Lemma~\ref{lem:ranktwo} becomes singular.
We therefore combine the periodic matching identity of
Theorem~\ref{thm:cayley}(1) with the signed current classification of
Proposition~\ref{thm:current-split} to analyze periodic eigenspaces on
$\mathcal T(H)=0$, and then apply coercivity and trigonometric bounds to the
polygonal modes.

\subsection{Common holomorphic solutions, reflection duality, and Krein signature}\label{subsec:periodic-common}

\begin{lemma}[Periodic matching and analytic continuation across the disk]\label{lem:periodic-common}
For the complex block \eqref{eq:complex-op} with $0<e<1$, one has
\begin{equation}\label{eq:cayley}
 \ker(M-I)\simeq\ker(A+B).
\end{equation}
If $A^2=0$ and $E=\ker A\cap\ker B$, then $\dim E\le\nu_1(M)\le 2\dim E$.
Moreover, whenever $\ru_1\neq 0$ and $\det H\neq 0$, every vector in $E$
defines a solution analytic in a neighborhood of the closed unit disk whose
position vanishes at $r$.
\end{lemma}
\begin{proof}
Identity \eqref{eq:cayley} and the bounds on $\nu_1(M)$ are immediate from
Theorem~\ref{thm:cayley}(1), and $I-\ii B/2$ acts as the identity on $E$. In local Frobenius
coordinates at either singularity, a fixed vector is annihilated by the
nilpotent logarithm because $\exp(\ii N_{\rm loc})-I$ equals $\ii N_{\rm loc}$ times an
invertible polynomial in $N_{\rm loc}$. At the origin, when $\ru_1\neq 0$, the
second row of \eqref{eq:generator} forces the exponent $-2$ coordinate to
vanish, leaving only nonnegative integer exponents with no logarithmic terms.
At $r$, invertibility of $H$ forces any logarithm-free fixed solution to be
analytic with position $w(r)=0$. A solution fixed by both based loops is
therefore single-valued on the punctured disk and extends holomorphically
across $0$ and $r$; since the next singularity lies at $r^{-1}<-1$, the
solution is analytic in a neighborhood of the closed unit disk.
\end{proof}

\begin{proposition}[Double periodic eigenspaces and Hermitian signature $(1,1)$]\label{prop:positive-current}
Suppose $a,b>0$, $c\neq 0$, $\mathcal T(H)=0$, and $\det H\neq 0$ in
\eqref{eq:complex-op}. Then at every periodic root with $0<e<1$, one has
$\nu_1(M)=2$, the multiplier $+1$ is semisimple of algebraic
multiplicity two, and the physical Hermitian form
$[\,\cdot,\cdot\,]_{\mathcal J}$ of Proposition~\ref{prop:current} has
signature $(1,1)$ on $\ker(M-I)$.
\end{proposition}
\begin{proof}
Since $\mathcal T(H)=0$, eliminating $c^2$ from \eqref{eq:T} gives the
identities
\begin{equation}\label{eq:LR}
\begin{aligned}
 L_a&=a-3b+3,\qquad R_a=3a-b-3,\qquad \det H=-\tfrac14 L_aR_a,\\
 3a(a-1)-c^2&=\tfrac34(a+b-1)R_a.
\end{aligned}
\end{equation}
Because $\det H\neq 0$ and $c\neq 0$, both $\ru_1$ and $\rv_2$ in
\eqref{eq:generator} are nonzero, while $\mathcal T(H)=0$ gives $A^2=0$; thus
$\rank A=\rank B=2$. Reversing the sign of the second component if needed makes
$c>0$ while preserving the physical Hermitian form. By
Lemma~\ref{lem:dual-transport}, the highest-exponent solution $w_+\in\ker A$
has strictly positive position at $r$, whereas every solution in $\ker B$
vanishes at $r$. Consequently $\ker A$ and $\ker B$ are distinct two-planes, so
at any periodic root their intersection $E$ is one-dimensional by
Lemma~\ref{lem:periodic-common}, which also implies $\nu_1(M)\le 2$.

Fix a nonzero $g_+\in E$; by Lemma~\ref{lem:periodic-common}, $g_+$ is
holomorphic past the closed unit disk with $g_+(r)=0$, so its Fourier expansion
on $|z|=1$ is supported in $k\ge 0$. Let $\Pi=\bigl(\begin{smallmatrix}0&1\\1&0\end{smallmatrix}\bigr)$
exchange the two components. Because $\Pi Z\Pi=-Z$ and $\q$ is even, the
involution $w(t)\mapsto\Pi w(-t)$ intertwines the periodic operators
$\mathcal L_{H,e}$ and $\mathcal L_{\Pi H\Pi,e}$. The swapped stiffness
$\Pi H\Pi$ satisfies the same hypotheses and shares the same periodic root $e$,
so it admits a nonzero holomorphic common solution $h_+$ vanishing at $r$. Its
reflection $g_-(t)=\Pi h_+(-t)$ then solves the original periodic equation with
Fourier support in $k\le 0$. If $g_+$ and $g_-$ were linearly dependent, $g_+$
would contain only the $k=0$ Fourier mode and hence be constant, contradicting
$g_+(r)=0$ and $g_+\not\equiv 0$. Thus $\{g_+,g_-\}$ is a basis of
$\ker(M-I)$ and $\nu_1(M)=2$.

Because $c\neq 0$, neither coordinate axis is an eigenspace of $H$, so the
$Z$-invariant unit-stiffness subspace $E_1$ in
Proposition~\ref{thm:current-split} vanishes. The one-sided Fourier supports
place $g_+$ and $g_-$ in the inner and outer periodic spaces, respectively,
whence Proposition~\ref{thm:current-split} yields both semisimplicity at $+1$
and the orthogonal signature relations
\[
 [g_+,g_+]_{\mathcal J}>0,\qquad
 [g_-,g_-]_{\mathcal J}<0,\qquad [g_+,g_-]_{\mathcal J}=0,
\]
establishing \hyperref[main:resonant]{Theorem~D}(1).
\end{proof}

\subsection{Periodic invertibility in the first mode of the square}\label{subsec:square}

\begin{proof}[Proof of {\hyperref[main:resonant]{Theorem~D}}\textup{(2)}]
For the regular $1+4$ polygon, \eqref{eq:sums} gives
$\gamma_4=(1+2\sqrt2)/4$ and $p_4=\sqrt2$. Substituting \eqref{eq:first-mode}
into \eqref{eq:T} shows that the unique positive root of $\mathcal T(H)=0$
solves $(112-48\sqrt2)m=113+8\sqrt2$, giving $m_*=(839+395\sqrt2)/496$. Since
$a,b>0$, $a\neq b$, and $c>0$ for all $m>0$, \hyperref[main:complex]{Theorem~C} leaves only $\lambda=1$
at $m=m_*$ to check at positive eccentricity.

At $m_*$, we have $A^2=0$ by \eqref{eq:K}, while
$a-3b+3=(9-4\sqrt2)/(4m_*+1+2\sqrt2)>0$ gives $\ru_1\neq 0$, and the box bounds
\eqref{eq:square-stiffness-bounds} yield $R_a=3a-b-3>0$ and hence $\det H<0$ by
\eqref{eq:LR}. If $0<e<1$ were a periodic root at $m_*$, then
Lemma~\ref{lem:periodic-common} would supply a nonzero common fixed solution
holomorphic past the closed unit disk and vanishing at $r$, whose transformed
trace under $\zeta=(z-r)/(1-rz)$ lies in the strictly positive Fourier sector.
Testing the periodic equation against $w$ gives $\mathfrak q_e[w]=0$,
contradicting the two-matrix positive-sector coercivity bound of
Lemma~\ref{lem:square-hardy} and Corollary~\ref{cor:square-box}. Since the
periodic operator is self-adjoint with compact resolvent, it is invertible for
every $0<e<1$; at $e=0$, Lemma~\ref{lem:square-circular} verifies invertibility
directly from the integer Fourier blocks.
\end{proof}

\subsection{Modal trichotomy for higher polygonal modes}\label{subsec:higher-periodic}

For each higher complex mode $2\le\ell<n/2$, introduce the auxiliary quantities
\begin{equation}\label{eq:exceptional-mass}
 u_\ell=\tfrac12(Q_{\ell-1}+Q_{\ell+1}),\qquad
 v_\ell=\gamma_n+2Q_\ell,\qquad
 \mu_{\mathcal T}=\frac{9v_\ell^2-12S_\ell^2}{12u_\ell+18v_\ell},\qquad
 m_{\mathcal T}=\mu_{\mathcal T}-\gamma_n.
\end{equation}
Throughout this subsection we abbreviate $u=u_\ell$ and $v=v_\ell$ while
retaining the subscript on $S_\ell$. By \eqref{eq:polygon-addition}, the
stiffness entries \eqref{eq:high-mode} read
\[
 a=\tfrac32+\frac{u-S_\ell}{\mu},\qquad
 b=\tfrac32+\frac{u+S_\ell}{\mu},\qquad
 c=\tfrac32\left(1-\frac{v}{\mu}\right),
\]
whence $\mu^2\mathcal T(H)=(12u+18v)\mu+12S_\ell^2-9v^2$. Because
$12u+18v>0$, \eqref{eq:exceptional-mass} defines the unique candidate root of
$\mathcal T(H(m))=0$.

\begin{lemma}[Stiffness bounds at an exceptional mass]\label{lem:exceptional-signs}
For every $2\le\ell<n/2$, one has $0<S_\ell\le u<v/2$ and
$v/4<\mu_{\mathcal T}<v/2$. In particular $c<0$ at $\mu=\mu_{\mathcal T}$.
At this mass all real shifts of the stiffness satisfy
\begin{equation}\label{eq:resonant-shift}
 \mu_{\mathcal T}^2\det(H(m_{\mathcal T})+hI)
 =(u+h\mu_{\mathcal T})^2+3h\mu_{\mathcal T}^2-4S_\ell^2,\qquad h\in\RR .
\end{equation}
\end{lemma}
\begin{proof}
Lemma~\ref{lem:polygon-sums} gives $S_\ell=(Q_{\ell+1}-Q_{\ell-1})/2>0$ and
$S_\ell\le u$. The second difference of the sums $Q_\ell$ gives
\[
 v-2u=\sum_{j=1}^{n-1}\frac{1-\cos(\ell\theta_j)}{2d_j}>0.
\]
In \eqref{eq:exceptional-mass}, the numerator is strictly greater than
$6v^2$ and the denominator is strictly less than $24v$. The numerator is
at most $9v^2$, and the denominator is strictly greater than $18v$. These
bounds give the claimed interval and the sign of $c$.
On $\mathcal T=0$, \eqref{eq:LR} gives
$\mu_{\mathcal T}^2\det H=u^2-4S_\ell^2$. Since
$\tr H=3+2u/\mu_{\mathcal T}$, expanding $\det(H+hI)=\det H+h\tr H+h^2$
gives \eqref{eq:resonant-shift}.
\end{proof}

\begin{lemma}[Absence of positive exceptional masses in the second mode]\label{lem:ell-two}
For every $n\ge 5$, the second complex mode satisfies $m_{\mathcal T}<0$.
\end{lemma}
\begin{proof}
Set $s_j=\sin(\pi j/n)$, $u_0=\pi/(2n)$, and
$x=(\sum_{j=1}^{n-1}s_j)/\gamma_n$. Since $\gamma_n=\frac14\sum_{j=1}^{n-1}s_j^{-1}$
and at least one $s_j<1$, we have $0<x<4$. The sine sum in the proof of
Lemma~\ref{lem:polygon-sums} gives $\sum s_j=\cot u_0$, so the first two
instances of \eqref{eq:cot-moments} yield
\[
 \frac{Q_2}{\gamma_n}=\frac{3-x}{2},\qquad
 \frac{Q_3}{\gamma_n}=4-\alpha x,\qquad
 \alpha=\frac12\left(3+\frac{\cot(3u_0)}{\cot u_0}\right).
\]
Since $0<u_0\le\pi/10$, strict convexity of $\tan$ on $[0,3u_0]$ gives
$0<\cot(3u_0)/\cot u_0<1/3$ and hence
$3/2<\alpha<5/3$. Because $Q_1=0$, we obtain $2u=2S_2=Q_3$ and
$v/\gamma_n=4-x$. Dividing the positive denominator of
\eqref{eq:exceptional-mass} by $\gamma_n$, the numerator of
$m_{\mathcal T}/\gamma_n$ reduces to
\[
 x\bigl[(30\alpha-54)+(9-3\alpha^2)x\bigr].
\]
Since $9-3\alpha^2>0$ and the bracket at $x=4$ equals
$-6(2\alpha-3)(\alpha-1)<0$, the numerator is strictly negative on $(0,4)$.
\end{proof}

\begin{lemma}[Divergence of the periodic Morse index near collision]\label{lem:collision}
If the least eigenvalue $\sigma_-$ of a fixed real symmetric matrix $H$ in
\eqref{eq:complex-op} satisfies $\sigma_-<-1/8$, then the complex periodic
Morse index of $\mathcal L_{H,e}$ tends to infinity as $e\uparrow 1$. In
particular, there exists a sequence of distinct eccentricities $e_j\uparrow 1$
at which $\mathcal L_{H,e_j}$ has a nontrivial periodic kernel.
\end{lemma}
\begin{proof}
Shift the angle to $x=t-\pi$ so that the limiting weight at $e=1$ becomes
$1-\cos x$. Let $\mathbf e_-\in\RR^2$ be a unit eigenvector of $H$ for
$\sigma_-$, and test the limiting quadratic form on
$w(x)=\sqrt{x}\,\psi(\log x)\mathbf e_-$, where $\psi\in C_c^\infty(-\infty,0)$
is real. Because $w$ is real-valued, the gyroscopic term
$-2\operatorname{Re}\langle w,Z(-\ii w')\rangle$ vanishes identically. Using
$1-\cos x\le x^2/2$, changing variables to $y=\log x$, and integrating the
cross term by parts bounds the limiting quadratic form from above by
\[
 \int_{\RR}\bigl(|\psi'(y)|^2+(1/4+2\sigma_-)|\psi(y)|^2\bigr)\dd y.
\]
Since $1/4+2\sigma_-<0$, dilating a fixed bump makes this integral strictly
negative, and disjoint translates toward $-\infty$ span an arbitrarily
high-dimensional negative subspace. For any finite collection of such test
functions, their supports in $x$ stay bounded away from $0$, so the weights
converge uniformly as $e\uparrow 1$ and the same subspace remains negative
definite for all $e$ sufficiently close to $1$.

On any compact subinterval of $[0,1)$, $\mathcal L_{H,e}$ is a norm-continuous
bounded potential perturbation of a fixed self-adjoint operator with compact
resolvent; its Morse index is therefore locally bounded and locally constant
away from periodic roots. Unbounded growth as $e\uparrow 1$ thus forces
periodic roots in every left neighborhood of $1$.
\end{proof}

\begin{lemma}[Comparison of the polygonal Fourier sums]\label{lem:high-sum-comparison}
For $2\le j\le\lfloor n/2\rfloor+1$, the sequence $Q_j/(j^2-1)$ is strictly
decreasing. For every $4\le\ell<n/2$,
\[
 0<Q_{\ell+1}<\frac{\ell+2}{\ell-2}Q_{\ell-1}\le 3Q_{\ell-1},
 \qquad 0<S_\ell<v/4,\qquad v>3\gamma_n.
\]
\end{lemma}
\begin{proof}
Write $M_n=\lfloor n/2\rfloor$ and $h_k=Q_{k+1}-Q_k$. The proof of
Lemma~\ref{lem:polygon-sums} gives strict decrease of $h_k/(2k+1)$ up to
$k=M_n$. Since $Q_1=0$, the quotient
$Q_j/(j^2-1)$ is the weighted mean of $h_k/(2k+1)$ for
$1\le k<j$ with weights $2k+1$. It too decreases strictly.
Comparing $j=\ell-1$ and $j=\ell+1$ gives the first inequality, using
Lemma~\ref{lem:polygon-sums} and the reflection $Q_{n-j}=Q_j$ at the
upper odd endpoint. Taking the ratio of the difference and the sum
gives $S_\ell/u<2/\ell$. Lemma~\ref{lem:exceptional-signs} then gives
$S_\ell<v/4$.

To bound $v$, use the formula $Q_3/\gamma_n=4-\alpha x$ in the proof of
Lemma~\ref{lem:ell-two}, where $\alpha<5/3$ and
$x=\cot(\pi/(2n))/\gamma_n$. The proof of Lemma~\ref{lem:polygon-sums}
shows that $\gamma_n/n$ is increasing. Here $n\ge9$. Pairing the
reflected terms in $\gamma_9$ and using $\csc t\ge 1/t+t/6$ gives
\[
 \frac{\gamma_n}{n}\ge\frac{\gamma_9}{9}
 \ge \frac1{2\pi}\left(\frac{25}{12}+\frac{5\pi^2}{243}\right)
 >\frac1{2\pi}\frac{20}{9}.
\]
The cosecant bound on $0<t\le\pi/2$ follows from
$(1+t^2/6)(1-t^2/6+t^4/120)\le 1$ and the Taylor upper bound for
$\sin t/t$. As $\cot(\pi/(2n))<2n/\pi$, we obtain $x<9/5$ and
$Q_3>\gamma_n$. The same lemma gives $h_k>0$ for $k<M_n$, so $Q_\ell\ge Q_3$
and $v>3\gamma_n$.
\end{proof}

\begin{proposition}[Full periodic positivity in all higher modes $\ell\ge4$]\label{prop:fixed-modes}
Let $4\le\ell<n/2$. For every $m>0$ with $H(m)>0$, there is
$\varepsilon(H(m))>0$ such that for all $0\le e<1$,
\[
 \mathfrak q_e[w]\ \ge\ \varepsilon(H(m))\int_0^{2\pi}\q^{-1}|w|^2\dd t
 \qquad \bigl(w\in H^1_{\rm per}([0,2\pi],\CC^2)\bigr).
\]
The candidate $m_{\mathcal T}$ is positive, and the periodic operator at
$m_{\mathcal T}$ is invertible for all $0\le e<1$.
\end{proposition}
\begin{proof}
By Lemmas~\ref{lem:exceptional-signs} and \ref{lem:high-sum-comparison},
\[
 \mu_{\mathcal T}>
 \frac{(9-12/16)v^2}{24v}
 =\frac{11}{32}v>\frac{33}{32}\gamma_n.
\]
In particular $m_{\mathcal T}>0$. Set $R_b=3b-a-3=-L_a$ with $L_a,R_a$
as in \eqref{eq:LR}. For every $m>0$, the expressions for $a,b$ give
\[
 R_a=\frac{3Q_{\ell-1}-Q_{\ell+1}}{\mu}>0,\qquad
 R_b=\frac{3Q_{\ell+1}-Q_{\ell-1}}{\mu}>0,
\]
where $Q_{\ell+1}>Q_{\ell-1}>0$. The inequality $Q_{\ell+1}<3Q_{\ell-1}$
also gives $2S_\ell<u$. At the exceptional mass, \eqref{eq:resonant-shift}
with $h=0$ now gives $\det H>0$, so $H>0$. Whenever
$H>0$, both diagonal entries exceed $1$ and the two strict inequalities
$R_a,R_b>0$ allow application of Lemma~\ref{lem:full-periodic-bound}.
This proves the bound and excludes the periodic kernel at $m_{\mathcal T}$.
\end{proof}

\begin{corollary}[Opposite mass slopes at a periodic double root]\label{cor:mass-phase}
Let $2\le\ell<n/2$, $m_{\mathcal T}>0$, and $\det H(m_{\mathcal T})\neq 0$,
and fix any $e_0\in(0,1)$ at which this complex block has the multiplier $+1$
at $m=m_{\mathcal T}$. In coordinates $\lambda=e^{2\pi\ii\theta}$, the local
real unit-multiplier set at fixed $e_0$ is the union of two real-analytic mass
curves $m_\pm(\theta)$ through $(0,m_{\mathcal T})$ whose derivatives at
$\theta=0$ have opposite signs. In particular, for all sufficiently small
nonzero real $m-m_{\mathcal T}$, the two multipliers near $+1$ are distinct,
simple, and on the unit circle.
\end{corollary}
\begin{proof}
Write $h=m-m_{\mathcal T}$ and let $M(h)$ be the $2\pi$-return in the physical
state coordinates $(w,-\ii w')$ of Proposition~\ref{prop:current}, which
preserves $\mathcal J$. By Lemma~\ref{lem:exceptional-signs}, $c<0$ at
$m_{\mathcal T}$ and $m_{\mathcal T}<m_s:=2Q_\ell$, while
\eqref{eq:high-mode}, \eqref{eq:polygon-addition}, and
Lemma~\ref{lem:polygon-sums} give diagonal entries
$3/2+Q_{\ell\mp 1}/(m_{\mathcal T}+\gamma_n)>0$.
Proposition~\ref{prop:positive-current} therefore shows that $+1$ is a
semisimple double eigenvalue of $M(0)$ whose eigenspace $K=\ker(M(0)-I)$ has
$\mathcal J$-signature $(1,1)$.

Using the physical $L^2(\dd t)$ inner product, note that the rescaled family
$(m+\gamma_n)\mathcal L_{H(m),e_0}$ is affine in $m$. For every nonzero
periodic solution $w$ at $m_{\mathcal T}$, comparison with the decoupled mass
$m_s$ in \eqref{eq:decoupled-pos} yields
\[
 (m_{\mathcal T}+\gamma_n)(m_s-m_{\mathcal T})
 \langle w,\partial_m\mathcal L_{H(m),e_0}|_{m_{\mathcal T}}w\rangle
  =(m_s+\gamma_n)\langle w,\mathcal L_{H(m_s),e_0}w\rangle>0,
\]
so the mass derivative defines a positive definite Hermitian form on $K$.

A small circle around $+1$ isolates its two multipliers counted with
multiplicity, and the associated Riesz projection yields a real-analytic frame
$F(h)$ for this complex two-dimensional spectral subspace. Writing $MF=FA$ and
$G=F^*\mathcal J F$, we have $A(0)=I$, signature $(1,1)$ for $G(0)$, and
$A^*GA=G$ for real $h$. The principal logarithm near $I$ defines
$L(h)=-\ii\log A(h)$, which satisfies $L(0)=0$ and $L^*G=GL$ because
$\log A^{-1}=-\log A$. For $\alpha\in\CC^2$, let $w_\alpha$ be the periodic
solution with initial state $F(0)\alpha$. Differentiating $MF=FA$ at $h=0$ and
pairing with $K$ eliminates $(M(0)-I)F'(0)$ because its range is
$\mathcal J$-orthogonal to $K$; since the only nonzero block of
$\mathcal J\partial_m\mathcal A$ in the physical first-order matrix
$\mathcal A$ is its upper-left block $-\ii q_{e_0}^{-1}\partial_m H$,
variation of constants at fixed initial state gives
\[
 \beta^*G(0)A'(0)\alpha=-\ii\beta^*B\alpha,\qquad
 \beta^*B\alpha=\int_0^{2\pi}q_{e_0}^{-1}w_\beta^*
                (\partial_mH)(m_{\mathcal T})w_\alpha\dd t,
\]
and hence $-G(0)L'(0)=B>0$.

Consequently, $L(h)/h$ extends real-analytically across $h=0$ to
$L'(0)=-G(0)^{-1}B$, which is similar to the Hermitian matrix
$-B^{1/2}G(0)^{-1}B^{1/2}$ and thus has two simple real eigenvalues
$\xi_-(0)<0<\xi_+(0)$. Because $L(h)/h$ is $G(h)$-self-adjoint for real $h$,
its characteristic polynomial has real coefficients for real $h$, so the
implicit function theorem produces real-analytic eigenvalue branches
$\xi_\pm(h)$ (real for real $h$). The two phases near zero are therefore
$\theta_\pm(h)=h\xi_\pm(h)/(2\pi)$, yielding two distinct simple unit
multipliers for small $h\neq 0$ that exhaust the spectral subspace; their
analytic inverses $m_\pm(\theta)$ have slopes $m_\pm'(0)=2\pi/\xi_\pm(0)$ of
opposite signs.
\end{proof}

\begin{proof}[Proof of {\hyperref[main:resonant]{Theorem~D}}\textup{(3)}]
Assume $m_{\mathcal T}>0$ and $\det H(m_{\mathcal T})\neq 0$. Since
Lemma~\ref{lem:exceptional-signs} gives $c<0$,
Proposition~\ref{prop:positive-current} and Corollary~\ref{cor:mass-phase}
show that every periodic root at $m_{\mathcal T}$ with $0<e<1$ is semisimple
of algebraic multiplicity two with Hermitian signature $(1,1)$ and admits two
real-analytic mass curves with opposite phase slopes at fixed eccentricity. At
$e=0$, the parameters $R_a=(2u-4S_\ell)/\mu_{\mathcal T}$ and
$R_b=(2u+4S_\ell)/\mu_{\mathcal T}$ satisfy $R_a>-4$, $R_b>0$ (by
Lemma~\ref{lem:exceptional-signs}), and $R_a\neq 0$ (since $\det H\neq 0$), so
Lemma~\ref{lem:square-circular} ensures invertibility of the circular operator
and analyticity in $e$ isolates all periodic roots at $m_{\mathcal T}$ in
$(0,1)$. If moreover $\det(H(m_{\mathcal T})+I/8)<0$, then $\tr H>0$ forces
$\det H<0$ and $\sigma_-<-1/8$, whence Lemma~\ref{lem:collision} produces a
sequence of distinct roots $e_j\uparrow 1$.

The three modal cases follow directly from Lemma~\ref{lem:ell-two}
($m_{\mathcal T}<0$ for $\ell=2$, $n\ge 5$), Lemma~\ref{lem:polygon-mass}
($m_{\mathcal T}>0$ and $\det(H(m_{\mathcal T})+I/8)<0$ for $\ell=3$,
$n\ge 100$), and Proposition~\ref{prop:fixed-modes} ($m_{\mathcal T}>0$ and
uniform periodic coercivity whenever $H(m)>0$, including at $m_{\mathcal T}$,
for $4\le\ell<n/2$). Finally, for every higher complex mode,
\eqref{eq:decoupled-pos} excludes periodic roots at the decoupling mass
$2Q_\ell$, whereas \hyperref[main:complex]{Theorem~C} forces algebraic
simplicity of $+1$ at every other positive mass distinct from
$m_{\mathcal T}$, completing the proof in cases~(a) and~(c).
\end{proof}

Table~\ref{tab:trichotomy} summarizes the multiplicity conclusions for the
specified blocks. Each multiplicity in the table is for one block in the
complex Hilbert-space convention of Section~\ref{subsec:setup}.

\begin{table}[ht]
\centering
\footnotesize
\setlength{\tabcolsep}{4pt}
\renewcommand{\arraystretch}{1.17}
\begin{tabular}{@{}>{\raggedright\arraybackslash}p{.195\textwidth}>{\raggedright\arraybackslash}p{.206\textwidth}>{\raggedright\arraybackslash}p{.175\textwidth}>{\raggedright\arraybackslash}p{.35\textwidth}@{}}
\toprule
Block & $\lambda\neq1$ & Exceptional mass & $\lambda=1$ \\
\midrule
First, $n=4$
& $\nu_\lambda\le1$
& $m_*>0$
& No root at $m_*$. Every root at other masses is algebraically simple. \\
First, $n\ge7$
& $\nu_\lambda\le1$; $-1$ simple if present
& None on $m>0$
& Every root is algebraically simple. \\
Second, $\ell=2$, $n\ge5$
& $\nu_\lambda\le1$ when coupled; $\nu_{-1}\le1$ at all masses
& $m_{\mathcal T}<0$
& Every root is algebraically simple. \\
Third, $\ell=3$, $n\ge100$
& $\nu_\lambda\le1$ when coupled; $\nu_{-1}\le1$ at all masses
& $m_{\mathcal T}>0$, $\sigma_-<-1/8$
& At $m_{\mathcal T}$, every root is semisimple double with signature $(1,1)$,
including a sequence $e_j\uparrow1$. Roots at other masses are simple. \\
$4\le\ell<n/2$
& $\nu_\lambda\le1$ when coupled; $\nu_{-1}\le1$ at all masses
& $m_{\mathcal T}>0$
& No roots where $H>0$, in particular at $m_{\mathcal T}$. Every root is
algebraically simple. \\
Real terminal, even $n$, $\ell=n/2$
& $\nu_{-1}\le1$
& ---
& Periodic dimension given by Proposition~\ref{prop:shear}. \\
\bottomrule
\end{tabular}
\caption{Floquet multiplicities for regular $(1+n)$-gons in the stated block
ranges, with $m>0$ and $0<e<1$. In a higher complex mode, coupled means
$m\neq2Q_\ell$.}
\label{tab:trichotomy}
\end{table}

\begin{remark}[Remaining finite polygon sizes]\label{rem:remaining-sizes}
In the first complex mode, the remaining small polygons $n\in\{3,5,6\}$ also
admit a unique positive root $m_*$ of $\mathcal T(H)=0$ (where
$\det(H(m_*)+I/8)>0$; direct evaluation shows that $n=6$ already falls inside
the full-domain coercivity region of Lemma~\ref{lem:full-periodic-bound},
whereas $n=3,5$ lie just outside the two-matrix box of
Lemma~\ref{lem:square-hardy}). Likewise, in the third mode the explicit bound
$n\ge 100$ of \hyperref[main:resonant]{Theorem~D}(3b) is chosen to keep the
hand-checkable proof of Lemma~\ref{lem:polygon-mass} elementary rather than to
locate the sharp onset of $\det(H(m_{\mathcal T})+I/8)<0$, which direct
evaluation places at $n=67$.
\end{remark}

\appendix

\section{Analytic normal forms and order-two resonance coefficients}\label{app:gauge}

Our two local analyses use slightly different gauge normalizations: in
Section~\ref{sec:real} we construct a real symplectic gauge via Hamiltonian
exponentials, whereas in Section~\ref{sec:complex} it is simplest to set every
resonant gauge entry to zero. The following lemma establishes convergence for
both constructions simultaneously.

\begin{lemma}[Convergence and uniqueness past the last resonance]\label{lem:gauge-rec}
Consider $\zeta G'=\mathcal A(\zeta)G-G\mathcal N(\zeta)$ with $\mathcal A$
analytic near $0$, $\mathcal A(0)=\Lambda$, $G(0)=I$, and
$\mathcal N(\zeta)=\Lambda+\zeta^2N_2+\zeta^4N_4$, where $\Lambda$ is
diagonalizable with eigenvalues $-2,0,0,2$ and $[\Lambda,N_j]=jN_j$. Once the
Taylor coefficients of $G$ and $N_2,N_4$ through order four are chosen
consistently, there exists a unique analytic continuation of $G$ to all higher
orders, depending uniformly on compact parameter sets on which $\mathcal A$,
the initial jet, and the nonresonant inverses remain bounded. For a regular
singular point with nilpotent residue $\Lambda$ and $\mathcal N=\Lambda$, the
same assertion holds from order zero.
\end{lemma}
\begin{proof}
Write $\mathcal A=\sum_{k\ge 0}\mathcal A_k\zeta^k$ and
$G=\sum_{j\ge 0}G_j\zeta^j$, setting $N_k=0$ for $k\notin\{2,4\}$. Equating
coefficients of $\zeta^j$ ($j\ge 1$) yields
\begin{equation}\label{eq:gauge-rec}
 (jI-\operatorname{ad}\Lambda)G_j=\Phi_j-N_j,\qquad
 \Phi_j=\sum_{k=1}^{j}\mathcal A_kG_{j-k}-\sum_{k=1}^{j-1}G_{j-k}N_k.
\end{equation}
In an eigenbasis of $\Lambda$, at a nonresonant entry $(p,q)$ we divide by
$j-\Lambda_{pp}+\Lambda_{qq}$, whereas at a resonant entry we set
$(N_j)_{pq}=(\Phi_j)_{pq}$ and specify $(G_j)_{pq}$ according to the chosen
normalization. For $j>4$, all denominators are nonzero; when $\Lambda$ is
nilpotent, $\operatorname{ad}\Lambda$ is nilpotent and
$jI-\operatorname{ad}\Lambda$ is invertible for every $j\ge 1$.

Vectorizing the matrix recursion as
$(jI-L_0)g_j=\sum_{k=1}^j L_k g_{j-k}$ with $\norm{L_k}\le CR^{-k}$ ($k\ge 1$),
we have $\norm{(jI-L_0)^{-1}}\le 2/j$ for all sufficiently large $j$. Fixing
$0<R_1<R$, if $\norm{g_i}\le C_0R_1^{-i}$ for all $i<j$, then
\[
 \norm{g_j}\le C_0R_1^{-j}\frac{2C}{j}\frac{R_1}{R-R_1}.
\]
Choosing $C_0$ large enough to bound the finitely many initial coefficients
gives convergence on $|\zeta|<R_1$, uniformly on compact parameter sets.

For Lemma~\ref{lem:real-local}, starting instead with successive Hamiltonian
exponentials $\exp(\zeta^j G_j)$ keeps the transformed system in
$\mathfrak{sp}(4,\RR)$ at every order, since $\operatorname{ad}\Lambda$
preserves this Lie algebra. The resulting formal infinite product satisfies
\eqref{eq:gauge-rec} for its own resonant matrices $N_2,N_4$; because
$jI-\operatorname{ad}\Lambda$ is invertible for $j>4$, its higher coefficients
coincide with the unique convergent solution determined by its order-four jet,
and the identity $G^TJ_4G=J_4$ persists from the formal series to the
convergent gauge.
\end{proof}

\begin{proof}[Evaluation of the order-two coefficients in Lemma~\ref{lem:resonance}]
In the state $Y=(w,\vang w)^T$, the coefficient matrix of $\vang Y$ is
$\bigl(\begin{smallmatrix}0&I\\-p_e(z)H&2Z\end{smallmatrix}\bigr)$ with
$p_e(z)=2z/(ez^2+2z+e)=2z/e-4z^2/e^2+O(z^3)$. In the exponent basis
$(-2,0,0,2)$, the change-of-basis matrix and transformed stiffness block are
\[
 S_{\rm eig}=\begin{pmatrix}0&1&0&1\\1&0&1&0\\0&0&0&2\\-2&0&0&0\end{pmatrix},\qquad
 V_H=S_{\rm eig}^{-1}\begin{pmatrix}0&0\\-H&0\end{pmatrix}S_{\rm eig},
\]
so that $\mathcal A_1=2V_H/e$ and $\mathcal A_2=-4V_H/e^2$. At order one, every
denominator in \eqref{eq:gauge-rec} is nonzero; writing
$d_{pq}=\Lambda_{pp}-\Lambda_{qq}$, the order-two resonant entries
($d_{pq}=2$) evaluate explicitly as
\[
 (N_2)_{pq}=\frac{4}{e^2}\left[-(V_H)_{pq}
 +\sum_{j=1}^{4}\frac{(V_H)_{pj}(V_H)_{jq}}{1-d_{jq}}\right].
\]
For direct verification, the four rows of $2V_H$ are
\[
 (b,c,b,c),\qquad (c,a,c,a),\qquad (-b,-c,-b,-c),\qquad (-c,-a,-c,-a).
\]
Substituting $(p,q)=(2,1),(3,1),(4,2),(4,3)$ yields $\ru_1,\ru_2,\rv_1,\rv_2$
in \eqref{eq:generator}. At order four, $N_4$ is a scalar multiple of $E_{41}$,
which is annihilated on both sides by $N_2$ and hence does not affect
\eqref{eq:K}.
\end{proof}

\section{Periodic coercivity and the circular spectrum}\label{app:hardy}

We use unnormalized circle measure $\dd t$ in the following estimates. With
$\mathsf D_t=-\ii\partial_t$, the quadratic form of the periodic
operator \eqref{eq:complex-op} on $H^1_{\rm per}([0,2\pi],\CC^2)$ is
\[
 \mathfrak q_e[w]=\int_0^{2\pi}
 \bigl(|\mathsf D_t w|^2-2\operatorname{Re}(w^*Z\mathsf D_t w)
       +\q^{-1}w^*Hw\bigr)\dd t.
\]
For $0<e<1$ put $\eta=\sqrt{1-e^2}$ and $\xi=-r$, so that $e\xi=1-\eta$ and
$e/\xi=1+\eta$. Passing to the M\"obius-adapted angle $\phi=\arg\zeta$ of
Section~\ref{sec:graph} ($\zeta=(z-r)/(1-rz)$, $\dd\phi/\dd t=\eta/\q=(1-e\cos\phi)/\eta$)
and setting $\widetilde w(\phi)=w(t)$ and $\mathsf D_\phi=-\ii\partial_\phi$
gives, on the full periodic form domain for every $e<1$,
\begin{equation}\label{eq:square-hardy-form}
 \begin{aligned}
 \eta\mathfrak q_e[w]&=\int_0^{2\pi}
 \bigl((1-e\cos\phi)|\mathsf D_\phi\widetilde w|^2
  -2\eta\operatorname{Re}(\widetilde w^*Z\mathsf D_\phi\widetilde w)
  +\widetilde w^*H\widetilde w\bigr)\dd\phi,\\
 \norm{\widetilde w}_{L^2}^2&=\eta\int_0^{2\pi}\q^{-1}|w|^2\dd t.
 \end{aligned}
\end{equation}

In the Fourier basis of $L^2(\dd\phi)$, the scalar operator
$S_e=\mathsf D_\phi(1-e\cos\phi)\mathsf D_\phi$
has diagonal $k^2$ and off-diagonal $-ek(k+1)/2$ between $k$ and $k+1$.
These edges vanish at $k=-1,0$, so the constant, strictly positive,
and strictly negative Fourier sectors of \eqref{eq:square-hardy-form}
are mutually orthogonal for the form. The constant has matrix $H$.
A sign change of the second component, if necessary, makes the
off-diagonal stiffness $-d$, where $d=|c|$. In the positive sector the
matrix $T$ has nonpositive off-diagonal entries. The negative sector is
the same positive-sector problem after reversing $k$ and exchanging $a,b$.

For a positive sequence $\Phi$ on the vertices $i=(j,k)$, with
$j\in\{1,2\}$ and $k\ge1$, and a finitely supported complex sequence $u$,
expansion of the squares gives the exact identity
\begin{equation}\label{eq:periodic-edge}
 \langle u,Tu\rangle=
 \sum_i\frac{(T\Phi)_i}{\Phi_i}|u_i|^2
 +\sum_{i<j}(-T_{ij})\Phi_i\Phi_j
 \left|\frac{u_i}{\Phi_i}-\frac{u_j}{\Phi_j}\right|^2 .
\end{equation}
The relation also holds with $T$ replaced by $T-\varepsilon I$.
In particular, $T\Phi\ge\varepsilon\Phi$ entrywise implies the lower
bound $\varepsilon$ for the positive-sector matrix. Finite Fourier sums
are dense in each sector in $H^1$, and \eqref{eq:square-hardy-form} is
continuous there, so the lower bound extends to the form domain.

\begin{lemma}[Coercivity on the full periodic domain]\label{lem:full-periodic-bound}
Suppose $a,b>1$, $ab>d^2$, $b<3(a-1)$ and $a<3(b-1)$. There is
$\varepsilon(H)>0$ such that for all $0\le e<1$ and
$w\in H^1_{\rm per}([0,2\pi],\CC^2)$,
\[
 \mathfrak q_e[w]\ge\varepsilon(H)\int_0^{2\pi}\q^{-1}|w|^2\dd t .
\]
\end{lemma}
\begin{proof}
First assume $0<e<1$ and $d>0$. For $k\ge1$ put $f_k=\xi^k$ and
$y=\eta k$. The tridiagonal recurrence, including the vanishing backward
edge at $k=1$, gives
\[
 (S_ef)_k=y f_k,\qquad (S_e(yf))_k=(2y-1)y f_k.
\]
Choose constants $d/b<C<a/d$ and $1/(a-1)<B<3C/d$. Both intervals are
nonempty since $ab>d^2$ and $b<3(a-1)$. For
$\Phi_k=((1+By)f_k,Cf_k)^T$, the components of $T\Phi$ divided by $f_k$
are
\[
 a-dC+\bigl((a-1)B-1\bigr)y,\qquad
 bC-d+(3C-dB)y .
\]
Thus $T\Phi\ge\varepsilon_+\Phi$ with
$\varepsilon_+=\min\{a-dC,a-1-B^{-1},b-d/C\}>0$.
Exchanging $a,b$ gives a positive bound $\varepsilon_-$ on the
negative Fourier sector. The constant sector has the positive bound
$\lambda_{\min}(H)$. Taking the minimum of these three bounds in
\eqref{eq:square-hardy-form} and \eqref{eq:periodic-edge} proves the
claim with a constant independent of $e$. For each fixed $w$, continuity
of the physical form and weight at $e=0$ gives the same bound there.
If $d=0$, the periodic scalar gauges $e^{\pm\ii t}$ and
\eqref{eq:picone} give the assertion directly with
$\varepsilon(H)=\min\{a-1,b-1\}>0$ for all $0\le e<1$.
\end{proof}

When $H$ is indefinite, Lemma~\ref{lem:full-periodic-bound} no longer applies;
on the positive Fourier sector, factoring out the zero at $r$ directly in
physical modes yields the following two-matrix criterion.

\begin{lemma}[Two-matrix coercivity on the positive Fourier sector]\label{lem:square-hardy}
For a constant Hermitian $2\times2$ stiffness $H$ with diagonal entries $a,b$
and upper off-diagonal entry $c$, suppose that both
$C_0=H+\diag(0,1/8)$ and $C_1=H+\diag(-1,17/8)$ are positive definite,
and set $\varepsilon=\min\{\lambda_{\min}(C_0),\lambda_{\min}(C_1)\}$.
For every $0\le e<1$ and every $w\in H^1_{\rm per}([0,2\pi],\CC^2)$ with
strictly positive Fourier support in the transformed angle $\phi$, one has
\begin{equation}\label{eq:hardy-test}
 \mathfrak q_e[w]\ge\varepsilon\int_0^{2\pi}\q^{-1}|w|^2\,\dd t .
\end{equation}
At $e=0$ use $r=0$ and $\phi=t$. The hypotheses hold, in particular,
if $a\ge19/10$, $b\ge3/2$, and $|c|\le7/4$.
\end{lemma}
\begin{proof}
Because the disk automorphism sends $r$ to $0$, strictly positive Fourier
support in $\phi$ allows the factorization $w(z)=(z-r)v(z)$ with
$v\in H^1_{\rm per}$ supported in the nonnegative physical modes: approximating
the transformed trace in $H^1$ by positive Fourier polynomials in $\zeta$,
composing with the automorphism, and dividing by $z-r$ yields quotients
holomorphic on the closed disk, and both operations are continuous in $H^1$ on
the unit circle for fixed $e<1$.

Write $v(z)=(2\pi)^{-1/2}\sum_{j\ge0}(-1)^j v_j z^j$ and put $\xi=-r$. Since
$\q=|z-r|^2/(1+\xi^2)$, the weighted squared norm of $w$ equals
$(1+\xi^2)\sum_j|v_j|^2$, so the potential term in
$\mathfrak q_e[w]/(1+\xi^2)$ acts simply as the constant matrix $H$. The
constant physical mode of $w$ has zero kinetic energy, while its mode $j+1$
($j\ge 0$) has coefficient $(-1)^j(v_j-\xi v_{j+1})$ and kinetic weights
$(j+1)(j-1)$ and $(j+1)(j+3)$ in the two components.

Let $P$ be the orthogonal projection onto $(1,-\xi,0,\ldots)/\sqrt{1+\xi^2}$
in the scalar sequence space. Because only the $j=0$ weight of the first
component is negative, its normalized kinetic form is bounded below by $-P$.
For any finitely supported scalar sequence $y$, since
$(k+1)-k\xi^2\ge 1$ for $0\le\xi<1$, telescoping together with the reverse
triangle inequality and Cauchy--Schwarz gives
\[
 \|y\|^2\le\sum_{j\ge0}(j+1)(|y_j|^2-\xi^2|y_{j+1}|^2)
 \le(1+\xi)\|y\|
      \left(\sum_{j\ge0}(j+1)^2|y_j-\xi y_{j+1}|^2\right)^{1/2}.
\]
Because $(1+\xi)^2(1+\xi^2)\le 8$, the weights $(j+1)^2$ yield a normalized
lower bound $I/8$, while the remaining weights $2(j+1)$ in $(j+1)(j+3)$
contribute at least $2P$ from $j=0$, so the second normalized kinetic form
dominates $I/8+2P$. As $P$ acts identically on both components and commutes
with $I\otimes H$, the combined form is bounded below by $C_1$ on
$\operatorname{ran}P$ and by $C_0$ on $\ker P$. Passing from finite Fourier
truncations of $v$ to $H^1$ proves \eqref{eq:hardy-test}. In the stated
stiffness box, the diagonal products in $C_0$ and $C_1$ are at least $247/80$
and $261/80$, respectively, both strictly exceeding $|c|^2\le 245/80$.
\end{proof}

\begin{corollary}[Exceptional square stiffness]\label{cor:square-box}
The stiffness at the exceptional first-mode mass $m_*$ for the regular
$1+4$ polygon satisfies the hypotheses of Lemma~\ref{lem:square-hardy}.
\end{corollary}
\begin{proof}
At $m_*$, \eqref{eq:first-mode} gives
\[
 a-1=\frac{2823+395\sqrt2}{1926+1286\sqrt2},\qquad
 b-1=\frac{839+891\sqrt2}{1926+1286\sqrt2}.
\]
These expressions are decreasing and increasing in $\sqrt2$,
respectively. Substituting $24/17<\sqrt2<10/7$ gives
\begin{equation}\label{eq:square-stiffness-bounds}
 \frac{19}{10}<a<\frac{191}{100},\qquad
 \frac{39}{25}<b<\frac{157}{100}.
\end{equation}
Since $\mathcal T(H)=0$ at $m_*$, equation \eqref{eq:T} gives
\[
 c^2=\tfrac34(a-b)^2+\tfrac32(a+b)-\tfrac94
 <\tfrac34\left(\frac7{20}\right)^2+\tfrac32\left(\frac{87}{25}\right)-\tfrac94
 =\frac{4899}{1600}<\left(\frac74\right)^2 .
\]
Since $b>39/25>3/2$, the hypotheses of Lemma~\ref{lem:square-hardy} hold.
\end{proof}

\begin{lemma}[Circular periodic spectrum on the resonant surface]\label{lem:square-circular}
Suppose $\mathcal T(H)=0$, set $R_a=3a-b-3$ and $R_b=3b-a-3$, and take $e=0$.
The periodic operator $\mathcal L_{H,0}$ is invertible if and only if neither
$R_a$ nor $R_b$ belongs to $\{-2j(j-1):j\in\mathbb Z\}$. At the exceptional
square mass $m_*$, $\mathcal L_{H,0}$ is invertible and has complex Morse index
two.
\end{lemma}
\begin{proof}
On the $k$th integer Fourier mode ($k\in\mathbb Z$), the block matrix
$\mathsf H_k$ has diagonal entries $k^2-2k+a$ and $k^2+2k+b$ and off-diagonal
entry $c$. Using \eqref{eq:LR}, its determinant factors as
\[
 \det\mathsf H_k=\tfrac14\bigl(R_a+2k(k-1)\bigr)\bigl(R_b+2k(k+1)\bigr),
\]
and the two quadratic sequences $\{2k(k-1)\}$ and $\{2k(k+1)\}$ have the same
range over $k\in\mathbb Z$. At $m_*$, \eqref{eq:square-stiffness-bounds} gives
$R_a>0$ and $-4<R_b<0$; thus $\det\mathsf H_k<0$ holds precisely for
$k\in\{-1,0\}$, whereas every other block has positive determinant and positive
trace $2k^2+a+b$. This excludes any zero eigenvalue and gives Morse index two.
\end{proof}

\section{Elementary rational bounds for finite polygons}\label{app:finite-n}

This appendix records the two finite-polygon rational verifications invoked in
Lemma~\ref{lem:polygon-sums} and \hyperref[main:resonant]{Theorem~D}(3b), so that every finite threshold
in the main text remains fully self-contained without computer assistance.

\begin{lemma}[Threshold comparisons for the first polygonal mode]\label{lem:finite-first-mode}
The trigonometric sums of \eqref{eq:sums} and \eqref{eq:first-mode} satisfy
$4\gamma_7+p_7>14$, $p_{11}<11$, and $p_{12}>12$.
\end{lemma}
\begin{proof}
Write $u_0=\pi/(2n)$, so that telescoping
$2\sin u_0\sin(2ju_0)=\cos((2j-1)u_0)-\cos((2j+1)u_0)$ yields
$p_n=4\gamma_n-\cot u_0$. We first note the elementary bounds
\[
 \frac1t+\frac t6+\frac{7t^3}{360}\le\csc t\le\frac1t+\frac t6+\frac{t^3}{21}
 \quad(0<t\le\pi/2),\qquad
 \frac1t-\frac t2\le\cot t\le\frac1t-\frac t3\quad(0<t\le\pi/4).
\]
Indeed, multiplying the lower $\csc t$ bound by $t-t^3/6+t^5/120\ge\sin t$
produces $1-t^6/540+7t^8/43200\le 1$, whereas multiplying the upper bound by
$t-t^3/6\le\sin t$ gives $1+t^4\bigl((1-t^2/6)/21-1/36\bigr)\ge 1$ because
$t^2<5/2$. Similarly, the lower $\cot t$ bound follows directly from
$\cos t\ge 1-t^2/2$ and $\sin t\le t$, while the upper bound holds because
$(1-t^2/3)\sin t-t\cos t$ vanishes at $0$ and has nonnegative derivative
$\frac13 t(\sin t-t\cos t)$.

Pairing $j$ with $n-j$ in $4\gamma_n$ and summing $1/j$, $j$, and $j^3$ over
the half-range yields
\begin{align*}
 4\gamma_7+p_7&\ge\frac{112}{3\pi}+\frac{25\pi}{42}+\frac{2\pi^3}{245}>14,\\
 p_{11}&\le\frac{847}{30\pi}+\frac{21\pi}{44}+\frac{150\pi^3}{9317}<11,\\
 p_{12}&\ge 1+\frac{154}{5\pi}+\frac{31\pi}{72}+\frac{35\pi^3}{6912}>12,
\end{align*}
where the midpoint term contributes $1$ to $4\gamma_{12}$ in the identity
$p_{12}=4\gamma_{12}-\cot(\pi/24)$. Substituting
$223/71<\pi<22/7$ in the appropriate direction in each term establishes all
three strict inequalities (and these bounds on $\pi$ themselves follow from
Machin's identity $\pi=16\arctan(1/5)-4\arctan(1/239)$ by truncating the
alternating series at $(6,1)$ and $(5,2)$ terms, respectively).
\end{proof}

\begin{lemma}[Uniform supercritical bound for the third mode ($n\ge 100$)]\label{lem:polygon-mass}
In the third complex mode ($\ell=3$), every integer $n\ge 100$ satisfies
$m_{\mathcal T}>0$ and $\det\bigl(H(m_{\mathcal T})+I/8\bigr)<0$.
\end{lemma}
\begin{proof}
Set $t_n=\pi/(2n)$, $\gamma=\gamma_n$, $x=\cot(t_n)/\gamma$, and
$r_j=\cot(jt_n)/\cot t_n$. For $n\ge 100$, pairing the first $49$ terms of
$4\gamma=\sum_{j=1}^{n-1}\csc(\pi j/n)$ with their reflections and using
$\csc s>1/s$ gives $\gamma>n H_{49}/(2\pi)$, where $H_j=\sum_{i=1}^j 1/i$.
Applying convexity of $s\mapsto 1/s$ to the blocks $8$--$15$, $16$--$31$, and
$32$--$49$ yields $H_{49}\ge H_7+16/23+32/47+4/9>22/5$, which together with
$\cot t_n<2n/\pi$ gives $0<x<10/11$.

Put $k=2u_3=Q_4+Q_2$, $\rho=v/k$ with $v=\gamma+2Q_3$, and
$\delta=2S_3/k=(Q_4-Q_2)/k$. By \eqref{eq:cot-moments},
\[
 K:=\frac{k}{\gamma}=9-Bx,\quad B=\frac{7+3r_3+r_5}{2},\qquad
 \rho=1+\frac{(1+r_3+r_5)x}{2K},\qquad \delta=1-\frac{3-x}{K}.
\]
Because $\tan$ is convex on $[0,5t_n]$, we have $0<r_j\le 1/j$ for $j=3,5$,
whence $B\le 41/10$ and $v/\gamma=9-(3+r_3)x>4$; Lemma~\ref{lem:exceptional-signs}
then gives $k<v$ and $\mu_{\mathcal T}>v/4>\gamma$, proving $m_{\mathcal T}>0$.
Moreover, both $x/(9-41x/10)$ and $(3-x)/(9-41x/10)$ increase on $[0,10/11]$,
so evaluating at $x=10/11$ yields
\[
 \rho\le\frac{197}{174}<\frac{17}{15},\qquad
 \delta\ge\frac{35}{58},\qquad \left(\frac{35}{58}\right)^2>\frac{29}{80}.
\]

To apply \eqref{eq:resonant-shift}, set $p=\mu_{\mathcal T}/k>0$. Since
$\delta^2>29/80>1/3$, equation \eqref{eq:exceptional-mass} implies
\[
 p=\frac{3\rho^2-\delta^2}{2(1+3\rho)}
   =\frac{\rho}{2}-\frac16+\frac{1-3\delta^2}{6(1+3\rho)}
   <\frac{17}{30}-\frac16=\frac25.
\]
Evaluating \eqref{eq:resonant-shift} at the collision shift $h=1/8$ therefore
gives
\[
 p^2\det(H(m_{\mathcal T})+I/8)
 =\left(\tfrac12+\tfrac p8\right)^2+\tfrac38p^2-\delta^2
 <\left(\tfrac12+\tfrac1{20}\right)^2+\tfrac3{50}-\delta^2
 =\tfrac{29}{80}-\delta^2<0.\qedhere
\]
\end{proof}

\subsection*{Declaration of generative AI and AI-assisted technologies in the manuscript preparation process}
During the preparation of this work, the authors used generative AI tools to
assist with language editing and the verification of algebraic calculations.
After using these tools, the authors reviewed and edited the manuscript as
needed and take full responsibility for the content of the publication.

\par\smallskip\noindent
\begingroup\small
\begin{minipage}[t]{0.46\textwidth}
\textsc{Xijun Hu}\\
School of Mathematics\\
Shandong University\\
Jinan 250100, P. R. China\\
\href{mailto:xjhu@sdu.edu.cn}{\texttt{xjhu@sdu.edu.cn}}
\end{minipage}\hfill
\begin{minipage}[t]{0.46\textwidth}
\textsc{Li Wu}\\
School of Mathematics\\
Shandong University\\
Jinan 250100, P. R. China\\
\href{mailto:vvvli@sdu.edu.cn}{\texttt{vvvli@sdu.edu.cn}}
\end{minipage}
\endgroup
\end{document}